\documentclass[10pt,a4paper]{article}
\usepackage[english]{babel}
\usepackage{amssymb}
\usepackage{amsfonts}
\usepackage{amsmath}
\usepackage{color}
\usepackage{amsthm}
\usepackage{graphicx}
\usepackage{appendix}
\usepackage{upgreek}
\usepackage{multirow}
\usepackage[a4paper]{geometry}
\usepackage{algorithm}
\usepackage{algorithmic}
\usepackage{tabularx}
\usepackage{stmaryrd} 
\usepackage{cleveref}

\newcommand{\be}{\begin{equation}}
\newcommand{\ee}{\end{equation}}
\newcommand{\bi}{\begin{itemize}}
\newcommand{\ei}{\end{itemize}}

\newtheorem{definition}{Definition}
\newtheorem{corollaire}{Corollary}
\newtheorem{lemme}{Lemma}

\newtheorem{theoreme}{Theorem}

\newcommand{\N}{\mathbb{N}}

\newcommand{\R}{\mathbb{R}}

\newcolumntype{C}{>{\centering}X}
\newcommand{\off}[1]{}

\newcounter{auxFootnote}
\title{Fast convergence of inertial dynamics with Hessian-driven damping under geometry assumptions}

\author {J.-F. Aujol\footnote{Univ. Bordeaux, Bordeaux INP, CNRS, IMB, UMR 5251, F-33400 Talence, France}
\and 
Ch. Dossal\footnote{IMT, Univ. Toulouse, INSA Toulouse, Toulouse, France}\setcounter{auxFootnote}{\value{footnote}}
\and
V.H. Ho\`ang\footnotemark[\value{auxFootnote}]
\and
H. Labarri\`ere\footnotemark[\value{auxFootnote}] 
\and 
A. Rondepierre\footnotemark[\value{auxFootnote}] \footnote{LAAS, Univ. Toulouse, CNRS, Toulouse, France}\\
\footnotesize{Jean-Francois.Aujol@math.u-bordeaux.fr,}\\
\footnotesize{\{Charles.Dossal,vhoang,Hippolyte.Labarriere,Aude.Rondepierre\}@insa-toulouse.fr}}

\begin{document}
\nocite{*}
\maketitle
\begin{abstract}
First-order optimization algorithms can be considered as a discretization of ordinary differential equations (ODEs) \cite{su2014differential}. In this perspective, studying the properties of the corresponding trajectories may lead to convergence results which can be transfered to the numerical scheme. In this paper we analyse the following ODE introduced by Attouch et al. in \cite{attouch2016fast}:
\begin{equation*}
    \forall t\geqslant t_0,~\ddot{x}(t)+\frac{\alpha}{t}\dot{x}(t)+\beta H_F(x(t))\dot{x}(t)+\nabla F(x(t))=0,
\end{equation*}
where $\alpha>0$, $\beta>0$ and $H_F$ denotes the Hessian of $F$. This ODE can be derived to build numerical schemes which do not require $F$ to be twice differentiable as shown in \cite{attouch2020first,attouch2021convergence}. We provide strong convergence results on the error $F(x(t))-F^*$ and integrability properties on $\|\nabla F(x(t))\|$ under some geometry assumptions on $F$ such as quadratic growth around the set of minimizers. In particular, we show that the decay rate of the error for a strongly convex function is $O(t^{-\alpha-\varepsilon})$ for any $\varepsilon>0$. These results are briefly illustrated at the end of the paper.
\end{abstract}


\section{Introduction}

This paper focuses on the study of the ODE defined by:
\begin{equation}
\label{eq:Hessian_ODE}
\forall t\geqslant t_0,~\ddot{x}(t)+\frac{\alpha}{t} \dot{x}(t)+\beta H_F(x(t)) \dot{x}(t)+ \nabla F(x(t))=0,
\tag{DIN-AVD}
\end{equation}
where $t_0>0$, $\alpha>0$, $\beta>0$, $x(t_0)\in\R^n$, $\dot x(t_0)=0$ and $F:\R^n\rightarrow\R$ is a convex and $C^2$ function whose gradient and Hessian are respectively denoted by $\nabla F$ and $H_F$. We consider that the function $F$ has a non empty set of minimizers $X^*$ and we denote $F^*=\min\limits_{x\in\R^n}F(x)$. The underlying motivation of this analysis lies in the minimization of the function $F$. 

In \cite{su2014differential}, Su et al. highlight the link between optimization methods and dynamical systems. In particular, this paper considers Nesterov's accelerated gradient method (NAGM) introduced in \cite{nesterov27method} as a discretization of the following ODE:
\begin{equation}
\label{eq:Nesterov_ODE}
\ddot{x}(t)+\frac{\alpha}{t} \dot{x}(t)+\nabla F(x(t))=0,
\tag{AVD}
\end{equation}
and shows that the trajectories defined by \eqref{eq:Nesterov_ODE} and NAGM have related properties. In fact, the authors of \cite{su2014differential} prove that $F(x(t))-F^*=\mathcal{O}(t^{-2})$ for any $\alpha\geqslant3$ and they provide a similar convergence rate for the iterates of NAGM. This continuous approach has been widely adopted in recent works leading to convergence results on optimization schemes such as NAGM \cite{attouch2018fast,jendoubi2015asymptotics,attouch2019rate,apidopoulos2021convergence,sebbouh2019nesterov,aujol:hal-03491527,Aujol2019} and the Heavy-ball method \cite{balti2016asymptotic,aujol2022convergence,aujol2021convergence}.

Alvarez et al. introduce in \cite{alvarez2002second} the Dynamical Inertial Newton-like system defined by:
\begin{equation}
\ddot{x}(t)+\alpha \dot{x}(t)+\beta H_F(x(t)) \dot{x}(t)+ \nabla F(x(t))=0,
\tag{DIN}
\end{equation}
which is a combination of the Newton dynamical system and the Heavy-ball with friction system. This ODE involves an Hessian-driven damping term which reduces the oscillations related to the heavy-ball system. 

In \cite{attouch2016fast}, Attouch et al. combine a similar Hessian-driven damping term to an asymptotic vanishing damping term resulting in \eqref{eq:Hessian_ODE}. The case $\beta=0$ corresponds to \eqref{eq:Nesterov_ODE} which is related to NAGM. In fact, \eqref{eq:Hessian_ODE} can be linked to the high-resolution ODE for NAGM introduced by Shi et al. in \cite{shi2021understanding}. The authors of \cite{attouch2016fast} prove that if $\beta>0$ and $\alpha\geqslant3$, the convergence rate of the trajectories is the same as \eqref{eq:Nesterov_ODE} and that \begin{equation}\int_{t_0}^{+\infty}t^2\|\nabla F(x(t))\|^2dt<+\infty.\end{equation} This additional result is significant as it ensures the fast convergence of the gradient and therefore a reduction of oscillations. This property of \eqref{eq:Hessian_ODE} is directly linked to the Hessian-driven damping and justifies the derivation of this ODE in order to define an associated numerical scheme. In \cite{attouch2020first}, Attouch et al. study a more general ODE:
\begin{equation}
\ddot{x}(t)+\frac{\alpha}{t}\dot{x}(t)+\beta(t)H_F(x(t))\dot x(t)+b(t)\nabla F(x(t))=0,
\label{eq:ODE_Hgen}
\end{equation}
and similar convergence results are given under some conditions on $\beta$ and $b$. The authors introduce numerical schemes derived from \eqref{eq:ODE_Hgen} which take advantage of the additional term such as the Inertial Gradient Algorithm with Hessian Damping (IGAHD):
\begin{equation}
\left\{
\begin{gathered}
x_k=y_{k-1}-s\nabla F(y_{k-1}),\\
y_k=x_k+\alpha_k(x_k-x_{k-1})-\beta\sqrt{s}(\nabla F(x_k)-\nabla F(x_{k-1}))-\frac{\beta\sqrt{s}}{k}\nabla F(x_{k-1}),
\end{gathered}\right.
\label{eq:IGAHD_intro}
\end{equation}
where $\alpha_k=\frac{k-1}{k+\alpha-1}$, $\alpha>0$, $\beta\geqslant0$ and $s>0$. It is proved in \cite{attouch2020first,attouch2021convergence} that if $\alpha\geqslant 3$, $s\leqslant \frac{1}{L}$ and $\beta\in (0,2\sqrt{s})$, then the sequence $(x_k)_{k\in\N}$ defined by \eqref{eq:IGAHD_intro} satisfies $F(x_k)-F^*=\mathcal{O}\left(k^{-2}\right)$ and
\begin{equation}
	\sum_{k\in\N} k^2\|\nabla F(x_k)\|^2<+\infty.
\end{equation}
Note that this algorithm only requires $F$ to be differentiable as the Hessian-driven damping is treated as the time derivative of the gradient term.

The convergence of the trajectories of \eqref{eq:Hessian_ODE} and \eqref{eq:Nesterov_ODE} were studied under additional geometry assumptions on $F$. Such hypotheses allow faster convergence rates to be found and provide a better understanding of the behaviour of trajectories. Attouch et al. prove in \cite[Theorem~3.1]{attouch2016fast} that if $F$ is $\mu$-strongly convex, $\alpha>3$ and $\beta>0$, then the trajectories of \eqref{eq:Hessian_ODE} satisfy:
\begin{equation*}
F(x(t))-F^*=\mathcal{O}\left(t^{-\frac{2\alpha}{3}}\right).
\end{equation*}
Similar convergence rates are given for the trajectories of \eqref{eq:Nesterov_ODE} in \cite[Theorem~4.2]{su2014differential} and \cite[Theorem~3.4]{attouch2018fast}.

In this work, we give an analysis of the trajectories of \eqref{eq:Hessian_ODE} on more general assumptions on the geometry of $F$ than strong convexity. We consider functions behaving as $\|x-x^*\|^\gamma$ for $\gamma\geqslant1$ where $x^*\in X^*$ and we provide convergence results on $F(x(t))-F^*$ and $\|\nabla F(x(t))\|$. This geometry assumption on $F$ was studied in \cite{Aujol2019} and \cite{aujol:hal-03491527} for the trajectories of \eqref{eq:Nesterov_ODE} associated to NAGM. This work aims to give a better understanding of the convergence of the trajectories of \eqref{eq:Hessian_ODE} in this setting. The next step would be to study the convergence of corresponding numerical schemes under these assumptions using a similar approach.

The main contributions of this work can be summarized as follows:
\begin{enumerate}
	\item Non asymptotic bound on $F(x(t))-F^*$ for the trajectories of \eqref{eq:Hessian_ODE} under the assumption that $F$ has a quadratic growth around its minimizers. The resulting convergence rate is asymptotically the fastest in the literature for \eqref{eq:Hessian_ODE} under this set of assumptions. 
	\item Strong integrability property on $\|\nabla F(x(t))\|$ and $F(x(t))-F^*$ under the same assumptions. Given the geometry of $F$, this improved integrability of the gradient has a direct influence on the convergence of the trajectories $F(x(t))-F^*$. We give an asymptotic convergence rate which is faster than $\mathcal{O}\left(t^{-\frac{2\alpha}{3}}\right)$ under a weaker assumption than strong convexity.
	\item Asymptotic convergence rate for $F(x(t))-F^*$ and improved integrability of the gradient for flat geometries of $F$, i.e functions behaving as $\|x-x^*\|^\gamma$ where $\gamma>2$. 
\end{enumerate}

\section{Preliminary: Geometry of convex functions}

Throughout this paper, we assume that $\R^n$ is equipped with the Euclidean scalar product $\langle \cdot , \cdot \rangle$ and the associated norm $\|\cdot\|$. As usual $B(x,r)$ denotes the open Euclidean ball with center $x\in\R^n$ and radius $r>0$. For any real subset $X\subset\R^n$, the Euclidean distance $d$ is defined as:
\begin{equation*}
\forall x\in\R^n,~d(x,X)=\inf\limits_{y\in X}\|x-y\|.
\end{equation*}

In this section, we introduce some geometry conditions that will be investigated later on. 
\begin{definition}
Let $F$ : $\mathbb{R}^{n} \rightarrow \mathbb{R}$ be a convex differentiable function having a non empty set of minimizers $X^*$. The function $F$ satisfies the assumption $\mathcal{G}_{\mu}^\gamma$ for some $\gamma\geqslant1$ and $\mu>0$ if for all $x \in \mathbb{R}^{n}$,
\begin{equation}
\frac{\mu}{2}d\left(x, X^{*}\right)^{\gamma} \leqslant F(x)-F^*.
\label{eq:H2}
\end{equation}
\end{definition}

The hypothesis $\mathcal{G}_\mu^\gamma$ is a growth condition on the function $F$ ensuring that it grows as fast as $\|x-x^*\|^\gamma$ around its set of minimizers. The case $\gamma=2$ corresponds to functions having a quadratic growth around their minimizers including strongly convex functions. As we consider convex functions, $\mathcal{G}_\mu^\gamma$ is directly related to the \L{}ojasiewicz inequality as stated in the following lemma. The proof is given in \cite[Proposition~3.2]{garrigos2017convergence}.
\begin{lemme}
Let $F:\R^n\rightarrow \R\cup \{+\infty\}$ be a convex differentiable function having a non empty set of minimizers $X^*$. Let $F^*=\inf F$. If $F$ satisfies $\mathcal{G}_\mu^\gamma$ for some $\gamma\geqslant2$ and $\mu>0$, then $F$ has a global \L{}ojasiewicz property with an exponent $1-\frac{1}{\gamma}$, i.e there exists $K>0$ such that:
\begin{equation}
\forall x\in \R^n,\quad K\left(F(x)-F^*\right)\leqslant \|\nabla F(x)\|^{\frac{\gamma}{\gamma-1}}.
\label{eq:Loja_gen}
\end{equation}
Specifically, if $F$ satisfies $\mathcal{G}_\mu^2$ for some $\mu>0$, then:
\begin{equation}
\forall x\in \R^n,\quad 2\mu\left(F(x)-F^*\right)\leqslant \|\nabla F(x)\|^2.
\label{eq:Lojasiewicz_grad}
\end{equation}
\label{lem:Loja}
\end{lemme}

The following assumption was used in \cite{cabot2009long,su2014differential,Aujol2019,apidopoulos2021convergence} and can be seen as a flatness condition.

\begin{definition}
Let $F$ : $\mathbb{R}^{n} \rightarrow \mathbb{R}$ be a convex differentiable function having a non empty set of minimizers $X^*$. The function $F$ satisfies the assumption $\mathcal{H}_\gamma$ for some $\gamma\geqslant1$ if for all $x^*\in X^*$,
\begin{equation}\forall x\in \R^n,\quad F(x)-F^* \leqslant \frac{1}{\gamma}\left\langle\nabla F(x), x-x^{*}\right\rangle.
\label{eq:H1}\end{equation}
The function $F$ satisfies the assumption $\mathcal{H}_\gamma^{loc}$ for some $\gamma\geqslant1$ if for all $x^*\in X^*$ there exists $\nu>0$ such that,
\begin{equation}
\forall x\in B(x^*,\nu),\quad F(x)-F^* \leqslant \frac{1}{\gamma}\left\langle\nabla F(x), x-x^{*}\right\rangle.
\label{eq:H1loc}
\end{equation}
\end{definition}

To have an intuition of the geometry of functions satisfying $\mathcal{H}_\gamma$, observe that the flatness property \eqref{eq:H1} implies that for any minimizer $x^*\in X^*$, there exists $M>0$ and $\nu>0$ such that:
\begin{equation}
\forall x\in B(x^*,\nu),~F(x)-F^*\leqslant M\|x-x^*\|^\gamma,
\end{equation}
see \cite[Lemma~2.2]{Aujol2019}. Therefore, this assumption ensures that it does not grow too fast around its set of minimizers. Note that $\mathcal{H}_1$ corresponds to convexity and it is therefore always satisfied in our setting.

\section{Convergence rates of \eqref{eq:Hessian_ODE} under geometry assumptions} 

In this section, we state fast convergence rates for \eqref{eq:Hessian_ODE} trajectories that can be achieved when $F$ satisfies geometry assumptions such as $\mathcal{G}^\gamma_\mu$ and $\mathcal{H}_\gamma$. The convergence results are given first for sharp geometries and then for flat geometries.

\subsection{Sharp geometries}

\subsubsection{Contributions}

We first consider $F$ as a convex $C^2$ function having a unique minimizer $x^*$ and satisfying $\mathcal{H}_\gamma$ and $\mathcal{G}^2_\mu$ for some $\gamma\geqslant1$. These assumptions gather convex functions having a quadratic growth around their minimizers and consequently strongly convex functions.  This set of hypotheses was considered in \cite{Aujol2019} and \cite{aujol:hal-03491527} to analyse the convergence of the trajectories of \eqref{eq:Nesterov_ODE}. In this setting, Aujol et al. show in \cite[Theorem~4.2]{Aujol2019} that if $\gamma\leqslant2$ and $\alpha>1+\frac{2}{\gamma}$, the solution of \eqref{eq:Nesterov_ODE} satisfies
\begin{equation}
F(x(t))-F^*=\mathcal{O}\left(t^{-\frac{2\alpha\gamma}{\gamma+2}}\right).
\label{eq:Optimal_1}
\end{equation}
In \cite[Theorem~5]{aujol:hal-03491527}, the authors give a non asymptotic convergence bound that recovers asymptotically \eqref{eq:Optimal_1}. Such time-finite rate allows to identify the dependency of the bound according to each parameter. In fact, Aujol et al. show that the bound on $F(x(t))-F^*$ is proportional to $\alpha^{\frac{2\alpha\gamma}{\gamma+2}}$. As a consequence, large values of $\alpha$ will not necessarily accelerate the convergence of the trajectories of \eqref{eq:Nesterov_ODE} despite ensuring a better asymptotical convergence rate.

We provide similar convergence results for \eqref{eq:Hessian_ODE} which are summarized in the following theorem. These claims are discussed below. The proof can be found in Section \ref{sec:proof_sharp1}.

\begin{theoreme}
Let $F: \mathbb{R}^{n} \rightarrow \mathbb{R}$ be a convex $C^2$ function having a unique minimizer $x^*$. Assume that $F$ satisfies $\mathcal{H}_\gamma$ and $\mathcal{G}_\mu^2$ for some $\gamma\geqslant1$ and $\mu>0$. Let $x$ be a solution of \eqref{eq:Hessian_ODE} for all $t\geqslant t_0$ where $t_0>0$, $\alpha>0$ and  $\beta>0$. Let $\lambda=\frac{2\alpha}{\gamma+2}$. Then, $\lambda<\alpha$ and we have that:
\begin{enumerate}
	\item if $\alpha>1+\frac{2}{\gamma}$, then for all $t\geqslant t_0+\beta(\alpha-\lambda)$,
\begin{equation}
	F(x(t))-F^*\leqslant \frac{K}{\left(t+\beta(\lambda-\alpha)\right)^{\frac{2\alpha\gamma}{\gamma+2}}}
	\label{eq:sharp_gen}
\end{equation}
	where $K$ depends on $t_0$, $\alpha$, $\beta$, $\gamma$ and $\mu$. In particular, if $t_0\leqslant \frac{\alpha r^*}{(\gamma+2)\sqrt{\mu}}$, then for all $t\geqslant \frac{\alpha r^{*}}{(\gamma+2) \sqrt{\mu}}+\beta(\alpha-\lambda)$, inequality \eqref{eq:sharp_gen} holds with
	\begin{equation}\noindent
	K= C_{1} e^{\frac{2 \gamma}{\gamma+2} C_{2}\left(\alpha-1-\frac{2}{\gamma}\right)}\left(1+\tfrac{\beta\gamma\sqrt{\mu}}{r^*}\right) E_{m}\left(t_{0}\right)\left(\tfrac{\alpha r^*}{(\gamma+2)\sqrt{\mu}}\right)^{\frac{2 \alpha \gamma}{\gamma+2}} ,
	\label{eq:sharp_fatbound}
	\end{equation}
	\item if $\alpha=1+\frac{2}{\gamma}$, then for all $t\geqslant t_0+\beta$,
	\begin{equation}
	F(x(t))-F^*\leqslant \left((t_0+\beta)^2+\frac{\lambda^2+\sqrt{\mu}}{\mu}\right)e^{\frac{\beta}{t_0}}\frac{E_m(t_0)}{t(t-\beta)},
\label{eq:sharp3}
	\end{equation}
	\item if $\alpha\geqslant 1+\frac{2}{\gamma}$, then
\begin{equation}
\int_{t_0}^{+\infty} u^{\frac{2\alpha\gamma}{\gamma+2}}\| \nabla F(x(u)) \|^{2} du < +\infty,
\end{equation}
\end{enumerate}
where:
\begin{itemize}\renewcommand{\labelitemi}{$\bullet$}
    \item $r^*$ is the unique positive real root of the polynomial :$$r \mapsto r^{3}-(1+C_0)r^{2}-2(1+\sqrt{2}) r-4,$$
    \item $C_0=\beta\dfrac{\sqrt{\mu}\gamma(\gamma \lambda-1)}{\gamma \lambda-2},$
    \item $C_{1}=\left(1+\frac{2}{r^{*}}\right)^{2}$,
    \item $C_{2}=\frac{1+C_0}{r^{*}}+\frac{1+\sqrt{2}}{r^{* 2}}+\frac{4}{3 r^{* 3}}$,
    \item $E_{m}:t\mapsto \left(1+\dfrac{\beta\alpha}{t}\right)\left(F(x(t))-F^{*}\right)+\dfrac{1}{2}\left\|\dot{x}(t)+\beta \nabla F(x(t))\right\|^{2}$.
\end{itemize}
\label{thm:sharp1}
\end{theoreme}

The first claim ensures that the trajectories of \eqref{eq:Hessian_ODE} have the same asymptotical convergence rate than the trajectories of \eqref{eq:Nesterov_ODE} if $\alpha>1+\frac{2}{\gamma}$, i.e $F(x(t))-F^*=\mathcal{O}\left(t^{-\frac{2\alpha\gamma}{\gamma+2}}\right)$. This rate is still valid if $\alpha=1+\frac{2}{\gamma}$ as stated in the third claim. Observe that as strongly convex functions satisfy $\mathcal{H}_1$ and $\mathcal{G}^2_\mu$, we give the same convergence rate as \cite{attouch2016fast} for this class of functions i.e $\mathcal{O}\left(t^{-\frac{2\alpha}{3}}\right)$. This rate is also achieved for weaker hypotheses such as the combination of convexity and $\mathcal{G}^2_\mu$.

In addition, we give a tight bound on $F(x(t))-F^*$ in \eqref{eq:sharp_fatbound} which highlights the influence of $\alpha$ and $\beta$ on the convergence of the trajectories. Note that if $\beta=0$, this bound is the same as that given in \cite[Theorem~5]{aujol:hal-03491527} for \eqref{eq:Nesterov_ODE}. As in \eqref{eq:Nesterov_ODE}, the upper bound is proportional to $\alpha^{\frac{2\alpha\gamma}{\gamma+2}}$ and consequently setting $\alpha$ too large may not be efficient. Moreover, by optimizing the bounds given in \eqref{eq:sharp_fatbound} and \eqref{eq:sharp3} according to $\beta$ for $\alpha\geqslant1+\frac{2}{\gamma}$, we get that the optimal value is $0$. However, setting $\beta>0$ ensures the fast convergence of the gradient as stated in the third statement. This property of \eqref{eq:Hessian_ODE} is not valid for \eqref{eq:Nesterov_ODE}, i.e for $\beta=0$. This integrability result is an improvement of 
\begin{equation}
\int_{t_0}^{+\infty} u^2\|\nabla F(x(u))\|^2 du<+\infty,
\end{equation}
proved by Attouch et al. in \cite{attouch2016fast} for convex functions.

Furthermore, as we consider that $F$ is convex and satisfies $\mathcal{G}^2_\mu$ for some $\mu>0$, Lemma \ref{lem:Loja} states that it has a \L{}ojasiewicz property with exponent $\frac{1}{2}$. More precisely, for all $x\in\R^n$,
\begin{equation}
2\mu(F(x)-F^*)\leqslant\|\nabla F(x)\|^2.
\end{equation}
Thus, the integrability property of the gradient automatically implies a similar property on the error as stated in the following theorem. The proof is given in Section \ref{sec:proof_sharp2}.

\begin{theoreme}
Let $F: \mathbb{R}^{n} \rightarrow \mathbb{R}$ be a convex $C^2$ function having a unique minimizer $x^*$ and satisfying $\mathcal{G}_\mu^2$. Let $x$ be a solution of \eqref{eq:Hessian_ODE} for all $t\geqslant t_0$ where $t_0>0$, $\alpha\geqslant3$ and $\beta>0$. Then, for any $\varepsilon\in (0,1)$,
\begin{equation}
F(x(t))-F^*=\mathcal{O}\left(t^{-\alpha+\varepsilon}\right),
\end{equation}
and
\begin{equation}
\int_{t_0}^{+\infty} u^{\alpha-\varepsilon}\left(F(x(u))-F^*\right) du < +\infty.
\label{eq:thm2}
\end{equation}
\label{thm:sharp2}
\end{theoreme}

The asymptotical rate given in Theorem \ref{thm:sharp2} is faster than that given for strongly convex functions in \cite{attouch2016fast} ($\mathcal{O}\left(t^{-\frac{2\alpha}{3}}\right)$). Moreover, the integrability result \eqref{eq:thm2} is significantly strong as it implies the following statements which are proved in Section \ref{sec:proof_sharp3}.

\begin{corollaire}
Let $F: \mathbb{R}^{n} \rightarrow \mathbb{R}$ be a convex $C^2$ function having a unique minimizer $x^*$ and satisfying $\mathcal{G}_\mu^2$. Let $x$ be a solution of \eqref{eq:Hessian_ODE} for all $t\geqslant t_0$ where $t_0>0$, $\alpha\geqslant3$ and $\beta>0$. Then, for any $\varepsilon\in (0,1)$, as $t\rightarrow+\infty$,
\begin{enumerate}
	\item
	\begin{equation}
	F(z(t))-F^*=o\left(t^{-\alpha-1+\varepsilon}\right),
	\end{equation}
	where $z:t\mapsto\frac{\int_{t/2}^tu^{\alpha-\varepsilon} x(u)du}{\int_{t/2}^tu^{\alpha-\varepsilon} du}$.
	\item 
	\begin{equation}
	\inf\limits_{u\in[t/2,t]}\left(F(x(u))-F^*\right)=o\left(t^{-\alpha-1+\varepsilon}\right).
	\end{equation}
	\item 
	\begin{equation}
		\liminf\limits_{t\rightarrow+\infty} ~t^{\alpha+1-\varepsilon}\log(t)(F(x(t))-F^*) =0,
	\end{equation}
	where $\liminf\limits_{t\rightarrow+\infty} f(t)=\lim\limits_{t\rightarrow+\infty}\left[\inf\limits_{\tau\geqslant t}f(\tau)\right]$.
\end{enumerate}
\label{cor:sharp2}
\end{corollaire}

We would like to point out that Theorem \ref{thm:sharp2} relies on Lemma \ref{lem:H1loc} which states that a convex $C^2$ function automatically satisfies the assumption $\mathcal{H}_{2-\delta}^{loc}$ for any $\delta\in (0,1]$. Note that the $C^2$ assumption is necessary to study \eqref{eq:Hessian_ODE} but the corresponding algorithms require only $F$ to be differentiable. Hence, such a strong result may not be valid in the discrete case without this $C^2$ assumption.

\subsubsection{Sketch of proof of Theorem \ref{thm:sharp1}}

As the proof of Theorem \ref{thm:sharp1} is technical, we give a brief overview of it in this section. The full proof can be found in Section \ref{sec:proof_sharp1}.

Theorem \ref{thm:sharp1} states three claims: the first two claims are upper bounds of $F(x(t))-F^*$ in the cases $\alpha>1+\frac{2}{\gamma}$ and $\alpha=1+\frac{2}{\gamma}$ and the third claim is an integrability property on $\|\nabla F(x(t))\|$ for all $\alpha\geqslant1+\frac{2}{\gamma}$.

The proof of each statement relies on the following Lyapunov energy:
\begin{equation}
\mathcal{E}(t) =\left(t^{2}+t \beta(\lambda-\alpha)\right)\left(F(x(t))-F^{*}\right)+\frac{1}{2}\left\|\lambda\left(x(t)-x^{*}\right)+t(\dot{x}(t)+\beta \nabla F(x(t)))\right\|^{2},
\end{equation}
where $\lambda=\frac{2\alpha}{\gamma+2}$. \\
The first step consists in showing that the energy $\mathcal{E}$ satisfies a differential inequality of the form:
\begin{equation}
\forall t\geqslant T, \quad \mathcal{E}^\prime (t)+\frac{p}{t+\beta(\lambda-\alpha)}\mathcal{E}(t)\leqslant\tilde\varphi(t+\beta(\lambda-\alpha))\mathcal{E}(t),
\end{equation}
for some $T\geqslant t_0$, $p\geqslant0$ and $\tilde\varphi:\R^+\rightarrow\R^+$. In practice, we get a stronger inequality involving an additional term: for all $t\geqslant T$,
\begin{equation}
\mathcal{E}^\prime (t)+\frac{p}{t+\beta(\lambda-\alpha)}\mathcal{E}(t)+\beta t(t+\beta(\lambda-\alpha)) \| \nabla F(x(t)) \|^{2}\leqslant\tilde\varphi(t+\beta(\lambda-\alpha))\mathcal{E}(t),
\label{eq:sharp_fourthclaim}
\end{equation}
Then by defining the energy $\mathcal{H}$ as
\begin{equation}
\mathcal{H}:t\mapsto\mathcal{E}(t)(t+\beta(\lambda-\alpha))^pe^{-\tilde\Phi(t+\beta(\lambda-\alpha))},
\end{equation}
where $\widetilde{\Phi}(t)=-\int_{t}^{+\infty} \widetilde{\varphi}(x) dx$, it follows that $\mathcal{H}$ is decreasing for all $t\geqslant T$.

This allows us to write that for all $t_1\geqslant T$ we have
\begin{equation}
\forall t\geqslant t_1,\quad F(x(t))-F^*\leqslant e^{-\tilde\Phi(t+\beta(\lambda-\alpha))}\frac{\mathcal{H}(t_1)}{(t+\beta(\lambda-\alpha))^{p+2}}.
\label{eq:sharp_template}
\end{equation}
From this inequality come the two first statements:
\begin{itemize}\renewcommand{\labelitemi}{$\bullet$}
	\item equation \eqref{eq:sharp_gen} given in the first claim follows from a trivial simplification of  \eqref{eq:sharp_template}. The bound specified in \eqref{eq:sharp_fatbound} relies on an optimization of $t_1$ in order to get the tightest control on $F(x(t))-F^*$. To do this, $t_1$ is chosen as the minimizer of $$t\mapsto(t+\beta(\lambda-\alpha))^pe^{-\tilde\Phi(t+\beta(\lambda-\alpha))}.$$
	Developing \eqref{eq:sharp_template} for this value of $t_1$ leads to the final bound.
	\item the second claim is a rewriting of \eqref{eq:sharp_template} in the case $\alpha=1+\frac{2}{\gamma}$ for $t_1=T$.
\end{itemize}

The proof of the third claim is based on the inequality \eqref{eq:sharp_fourthclaim} and follows the same approach. By defining $\mathcal{G}$ as
\begin{equation}
\mathcal{G}:t\mapsto \mathcal{H}(t)+\beta\int_{T}^tu(u+\beta(\lambda-\alpha))^{p+1}e^{-\widetilde{\Phi}(u+\beta(\lambda-\alpha))}\|\nabla F(x(u))\|^2 du,
\end{equation}
the inequality \eqref{eq:sharp_fourthclaim} implies that $\mathcal{G}$ is decreasing and as $\mathcal{H}$ is a positive function,
\begin{equation}
\forall t\geqslant T,\quad \beta\int_{T}^tu(u+\beta(\lambda-\alpha))^{p+1}e^{-\widetilde{\Phi}(u+\beta(\lambda-\alpha))}\|\nabla F(x(u))\|^2 du\leqslant\mathcal{G}(T).
\end{equation}
Simple calculations lead to the final result.

\subsection{Flat geometries}

\subsubsection{Contributions}

We now focus on functions satisfying $\mathcal{H}_{\gamma_1}$ and $\mathcal{G}^{\gamma_2}_\mu$ where $\gamma_1\geqslant\gamma_2>2$ and $\mu>0$. These functions are said to have a flat geometry because they behave as $\|x-x^*\|^{\gamma_1}$ around their set of minimizers and $\gamma_1>2$. This geometry assumption was investigated in \cite{Aujol2019} for \eqref{eq:Nesterov_ODE}. The authors prove that if $\alpha\geqslant \frac{\gamma_1+2}{\gamma_1-2}$, then
\begin{equation}
F(x(t))-F^*=\mathcal{O}\left(t^{-\frac{2\gamma_2}{\gamma_2-2}}\right).
\label{eq:rate_flat_N}
\end{equation}
We provide a similar result for \eqref{eq:Hessian_ODE} in the following theorem which is proved in Section \ref{sec:proof_flat1}. We also give an integrability result on $\nabla F$ related to the Hessian driven damping term.

\begin{theoreme}
Let $F: \mathbb{R}^{n} \rightarrow \mathbb{R}$ be a convex $C^2$ function having a unique minimizer $x^*$. Assume that $F$ satisfies $\mathcal{H}_{\gamma_1}$ and $\mathcal{G}_\mu^{\gamma_2}$ for some $\gamma_1>2$, $\gamma_2>2$ such that $\gamma_{1}\geqslant\gamma_{2}$ and $\mu>0$. Let $x$ be a solution of \eqref{eq:Hessian_ODE} for all $t\geqslant t_0$ where $t_0>0$, $\alpha\geqslant\frac{\gamma_{1}+2}{\gamma_{1}-2}$ and $\beta>0$. Then as $t\rightarrow+\infty$,
\begin{equation}
F(x(t))-F^{*}=\mathcal{O}\left(t^{-\frac{2 \gamma_{1}}{\gamma_{1}-2}}\right).
\label{eq:rate_flat}\end{equation}
Moreover,
\begin{equation}
\int_{t_0}^{+\infty} u^{\frac{2\gamma_1}{\gamma_1-2}}\|\nabla F(x(u))\|^2du<+\infty.
\label{eq:integ_flat}
\end{equation}
\label{thm:flat}
\end{theoreme}
The asymptotical convergence rate \eqref{eq:rate_flat} is slightly slower than \eqref{eq:rate_flat_N} given by Aujol et al. for \eqref{eq:Nesterov_ODE} as $\gamma_1\geqslant\gamma_2$. However, we give an additional result on the integrability of the gradient which ensures a reduction of oscillations. As specified for sharp geometries, the assumption $\mathcal{G}^{\gamma_2}_\mu$ is equivalent to a \L{}ojasiewicz property with exponent $1-\frac{1}{\gamma_2}$. Consequently, we get that 
\begin{equation}
\int_{t_0}^{+\infty}u^{\frac{2\gamma_1}{\gamma_1-2}}\left(F(x(u))-F^*\right)^{\frac{\gamma_2-1}{2\gamma_2}}du<+\infty.
\end{equation}
This statement may lead to improved convergence rates according to the value of $\gamma_1$ and $\gamma_2$ as stated in the following corollary which is proved in Section \ref{sec:proof_flat_cor}.

\begin{corollaire}
Let $F: \mathbb{R}^{n} \rightarrow \mathbb{R}$ be a convex $C^2$ function having a unique minimizer $x^*$. Assume that $F$ satisfies $\mathcal{H}_{\gamma_1}$ and $\mathcal{G}_\mu^{\gamma_2}$ for some $\gamma_1>2$, $\gamma_2>2$ such that $\gamma_{1}\geqslant\gamma_{2}$ and $\mu>0$. Let $x$ be a solution of \eqref{eq:Hessian_ODE} for all $t\geqslant t_0$ where $t_0>0$, $\alpha\geqslant\frac{\gamma_{1}+2}{\gamma_{1}-2}$ and $\beta>0$. Then as $t\rightarrow+\infty$,
\begin{equation}
\inf\limits_{u\in[t/2,t]}F(x(u))-F^*=o\left(t^{-\frac{(3\gamma_1-2)\gamma_2}{2(\gamma_1-2)(\gamma_2-1)}}\right).
\end{equation}
\label{cor:flat1}
\end{corollaire}
Note that if $\gamma_1> 2$ and $\gamma_2\in\left(2,\frac{4\gamma_1}{\gamma_1+2}\right)$, then $\frac{(3\gamma_1-2)\gamma_2}{2(\gamma_1-2)(\gamma_2-1)}>\frac{2\gamma_1}{\gamma_1-2}$. Consequently, for this set of parameters the asymptotical rate of $\inf\limits_{u\in[t/2,t]}F(x(u))-F^*$ given in Corollary \ref{cor:flat1} is faster than the rate of $F(x(t))-F^*$ given in Theorem \ref{thm:flat}. 

\subsubsection{Sketch of proof of Theorem \ref{thm:flat}}

In this section, we give an outline of the proof of Theorem \ref{thm:flat} which is given in Section \ref{sec:proof_flat1}. The proof relies on the analysis of the Lyapunov energies $\mathcal{E}$ and $\mathcal{H}$ defined as:
\begin{equation}
\begin{gathered}
\mathcal{E}(t) =\left(t^{2}+t \beta(\lambda-\alpha)\right)\left(F(x(t))-F^{*}\right)+\frac{\xi}{2}\left\|x(t)-x^{*}\right\|^{2}\\+\frac{1}{2}\left\|\lambda\left(x(t)-x^{*}\right)+t(\dot{x}(t)+\beta \nabla F(x(t)))\right\|^{2},
\end{gathered}
\end{equation}
\begin{equation}
\mathcal{H}(t)=t^p\mathcal{E}(t),
\end{equation}
where $x^*$ is the unique minimizer of $F$, $\lambda\in\R$, $\xi\in\R$ and $p>0$.

The first step is to show that for a well-chosen set of parameters $\left(\lambda,\xi,p\right)$, the following inequality holds:
\begin{equation}
\forall t\geqslant t_1,~\mathcal{H}^\prime(t)+\beta t^{p+1}(t+\beta(\lambda-\alpha))\|\nabla F(x(t))\|^2\leqslant t^p\beta C_1\left(F(x(t))-F^*\right),
\label{eq:flat_sop_ine}
\end{equation}
where $t_1>0$ and $C_1>0$. In this case, the right term is not zero, implying that we cannot directly deduce that $\mathcal{H}$ is decreasing.\\
Therefore, the second step of the proof consists in investigating the function $\mathcal{G}$ defined by:
\begin{equation}
\mathcal{G}:t\mapsto \mathcal{H}(t)-\beta C_1 \int_{t_1}^t u^{p}(F(x(u))-F^*)du,
\end{equation}
where $t_1>T$ and $T$ is a well-chosen parameter. The objective is to use the decreasing nature of $\mathcal{G}$ to show that:
\begin{equation}
F(x(t))-F^*=\mathcal{O}\left(\frac{1}{t^{p+2}}\right).
\end{equation}
For this purpose, we show that the function $v$ defined by:
\begin{equation}
v(t)= t(t+\beta(\lambda-\alpha))^{p+1}(F(x(t))-F^*),
\end{equation}
is bounded by using the assumptions of the theorem.\\
The third step is to prove the second statement namely:
\begin{equation}
\int_{t_0}^{+\infty} u^{\frac{2\gamma_1}{\gamma_1-2}}\|\nabla F(x(u))\|^2du<+\infty.
\label{eq:flat_sop_grad}
\end{equation}
This is done by introducing the function $\mathcal{F}$ defined as follows:
\begin{equation}
\mathcal{F}:t\mapsto \mathcal{H}(t)-\beta C_1 \int_{t_1}^t u^{p-1}a(u)du+\beta \int_{t_1}^t u^{p+1}(u+\beta(\lambda-\alpha))\|\nabla F(x(u))\|^2du.
\end{equation}
Equation \eqref{eq:flat_sop_grad} is obtained by combining the decreasing nature of $\mathcal{F}$ and the boundedness of $v$.

\section{Numerical experiments}

In this section, we illustrate the fast convergence rates obtained theoretically for \eqref{eq:Hessian_ODE} with numerical experiments. We consider the following least-squares problem:
\begin{equation}
\min\limits_{x\in\R^n}F(x):=\|Ax-b\|^2,
\label{eq:Least-squares}
\end{equation}
where $A\in\mathcal{M}_{n\times n}\left(\R\right)$ and $b\in\R^n$. The function $F$ is convex, $C^2$ and satisfies $\mathcal{G}^2_\mu$ for some $\mu>0$. We apply the Inertial Gradient Algorithm with Hessian Damping (IGAHD) which was introduced by Attouch et al. in \cite{attouch2020first}:
\begin{equation}
\left\{
\begin{gathered}
x_k=y_{k-1}-s\nabla F(y_{k-1}),\\
y_k=x_k+\alpha_k(x_k-x_{k-1})-\beta\sqrt{s}(\nabla F(x_k)-\nabla F(x_{k-1}))-\frac{\beta\sqrt{s}}{k}\nabla F(x_{k-1}),
\end{gathered}\right.
\label{eq:IGAHD2}
\end{equation}
where $\alpha_k=\frac{k-1}{k+\alpha-1}$, $\alpha>0$, $\beta\geqslant0$ and $s>0$. Observe that the case $\beta=0$ corresponds to the Nesterov's accelerated gradient method. This numerical scheme is derived from the following ODE:
\begin{equation}
\label{eq:Hessian_modif}
\forall t\geqslant t_0,~\ddot{x}(t)+\frac{\alpha}{t} \dot{x}(t)+\beta H_F(x(t)) \dot{x}(t)+ \left(1+\frac{\beta}{t}\right)\nabla F(x(t))=0,
\end{equation}
which is a slightly modified version of \eqref{eq:Hessian_ODE}. The additional vanishing coefficient in front of the gradient keeps the structure of the dynamic while facilitating the computational aspects.

We do not provide any convergence results on IGAHD in this paper but we want to emphasize that the convergence of the iterates of IGAHD is related to the convergence rates obtained for \eqref{eq:Hessian_ODE}. We refer the reader to \cite{attouch2020first,attouch2021convergence} for a detailed analysis of this method. We recall that Attouch et al. prove in \cite{attouch2020first} that if $\alpha\geqslant3$, $s\leqslant\frac{1}{L}$ and $\beta\leqslant 2\sqrt{s}$, then the sequence $(x_k)_{k\in\N}$ defined in \eqref{eq:IGAHD2} satisfies:
\begin{equation}
F(x_k)-F^*=\mathcal{O}\left(k^{-2}\right).
\end{equation}

\begin{figure}[H]
\centering
\includegraphics[width=1\textwidth]{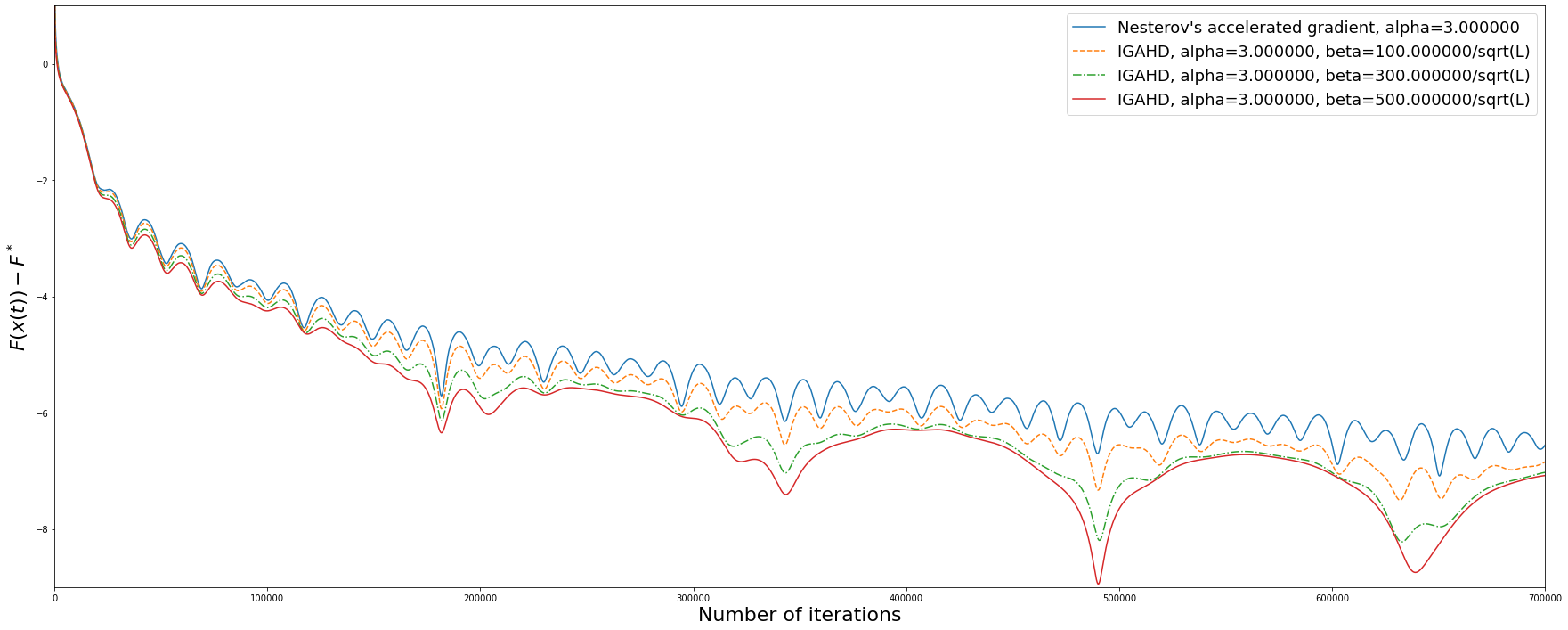}
\caption{Convergence rate of IGAHD for several values of $\beta$ for a least squares problem ($N=500$).}
\label{fig:multi_beta_LS}
\end{figure}

We compare the convergence of the iterates of IGAHD for several values of $\beta$ with the iterates of Nesterov's accelerated gradient method ($\beta=0$) to observe the influence of the Hessian-driven damping. Figure \ref{fig:multi_beta_LS} shows that the additional Hessian related term has a significant impact on the oscillations of the iterates. Indeed, this pathological behavior is reduced as $\beta$ grows. This can be related to the fast convergence of the gradient demonstrated in Theorem \ref{thm:sharp1}. 
Note that these experiments were made for large values of $\beta$ ($\beta\geqslant \frac{100}{\sqrt{L}}$) and consequently the convergence results given in \cite{attouch2020first,attouch2021convergence} do not hold in this context. Moreover, $\beta$ cannot be chosen too large as the iterates may not converge. There exists a critical value $\tilde\beta$ such that the algorithm does not converge for all $\beta\geqslant \tilde \beta$ and this value vary according to the geometry of $F$. However, no theoretical result on $\tilde\beta$ has been proved.

\section{Proofs}
\subsection{Proof of Theorem \ref{thm:sharp1}}

\label{sec:proof_sharp1}

Let $\alpha\geqslant 1+\frac{2}{\gamma}$ and $\lambda=\frac{2 \alpha}{\gamma+2}$. We consider the following Lyapunov function:
\begin{equation*}
\mathcal{E}(t) =\left(t^{2}+t \beta(\lambda-\alpha)\right)\left(F(x(t))-F^{*}\right)+\frac{1}{2}\left\|\lambda\left(x(t)-x^{*}\right)+t(\dot{x}(t)+\beta \nabla F(x(t)))\right\|^{2}.
\end{equation*}

\begin{lemme}
For all $t\geqslant t_0$:
\begin{equation}
\begin{aligned}\mathcal{E}^{\prime}(t) &=(2-\gamma \lambda)(t+\beta(\lambda-\alpha))\left(F(x(t))-F^{*}\right)-\beta(\lambda-\alpha)\left(F(x(t))-F^{*}\right)\\ 
&-\lambda(t+\beta (\lambda -\alpha))\left[-\gamma\left(F(x(t))-F^{*}\right)+\left\langle\nabla F(x(t)), x(t)-x^{*}\right\rangle\right]\\ 
&+\frac{\lambda(\lambda+1-\alpha)}{t}\left\langle x(t)-x^{*}, t(\dot{x}(t)+\beta \nabla F(x(t)))\right\rangle \\ 
& +\frac{\lambda+1-\alpha}{t}\|t(\dot{x}(t)+\beta \nabla F(x(t)))\|^{2}- \beta t(t+\beta(\lambda+\alpha)) \| \nabla F(x(t))\|^{2}.
\end{aligned}
\end{equation}
The proof of this lemma is given in Section \ref{subsection:Proof for Lemma 1}.
\label{lem:sharp1}
\end{lemme}

By using Lemma \ref{lem:sharp1} and the assumptions on $F$ we get the following result.
\begin{lemme}
Let $t\geqslant\max\{t_0,\beta(\alpha-\lambda)\}$. Then,
\begin{equation}
\begin{aligned}
\mathcal{E}^{\prime}(t)+\frac{\gamma \lambda-2}{t} \mathcal{E}(t) & \leqslant K(\alpha)\left( \;\dfrac{\lambda}{t}\left\|x(t)-x^{*}\right\|^{2}+\left\langle x(t)-x^{*}, \dot{x}(t)+\beta \nabla F(x(t))\right\rangle \right)
\\
&+\beta(\alpha-\lambda)\left(F(x(t))-F^{*}\right)-\beta t(t+\beta(\lambda-\alpha)) \| \nabla F(x(t)) \|^{2},
\end{aligned}
\label{eq:sharp1}
\end{equation}
where  $K(\alpha)=\dfrac{2 \alpha \gamma}{(\gamma+2)^{2}}\left(\alpha-1-\frac{2}{\gamma}\right)=\dfrac{(\gamma \lambda-2)\lambda}{2}$. In particular, if $\alpha=1+\frac{2}{\gamma}$, then
\begin{equation}
\mathcal{E}^{\prime}(t)+\beta t(t+\beta(\lambda-\alpha))\|\nabla F(x(t))\|^2\leqslant \beta(\alpha-\lambda)(F(x(t))-F^*).
\end{equation}
The proof of this lemma is given in Section \ref{subsection:Proof for Lemma 3}
\label{lem:sharp3}
\end{lemme}
\noindent\textbf{Case $\alpha>1+\frac{2}{\gamma}$ (Proof of statements 1 and 3).} The inequality \eqref{eq:sharp1} can be rewritten in the following way.
\begin{lemme}
Let $t>\max\{t_0,\beta(\alpha-\lambda)\}$. Then, if $\alpha>1+\frac{2}{\gamma}$,
\begin{equation}
\mathcal{E}^{\prime}(t)+\frac{\gamma \lambda-2}{t+\beta(\lambda-\alpha)} \mathcal{E}(t) +\beta t(t+\beta(\lambda-\alpha)) \| \nabla F(x(t)) \|^{2}\leqslant \widetilde{\varphi}\left(t+\beta(\lambda-\alpha)\right) \mathcal{E}(t),
\end{equation}
where 
$$
\widetilde{\varphi}:t\mapsto\frac{K(\alpha)}{\mu t^{2}}\left(\sqrt{\mu}(1+C_0)+\frac{2 \alpha}{(\gamma+2) t}(1+\sqrt{2})+\frac{4 \alpha^{2}}{(\gamma+2)^{2} \sqrt{\mu} t^{2}}\right),
$$
and $C_0=\beta\dfrac{\sqrt{\mu}\gamma(\gamma \lambda-1)}{\gamma \lambda-2}.$\\
The proof of this lemma is given in Section \ref{subsection:Proof for Lemma 4}
\label{lem:sharp4}
\end{lemme}
Let $\mathcal{H}$ be defined as follows:
$$
\mathcal{H}: t \mapsto \mathcal{E}(t) (t+\beta(\lambda-\alpha))^{\gamma \lambda-2} e^{-\widetilde{\Phi}(t+\beta(\lambda-\alpha))},$$
where $\widetilde{\Phi}(t)=-\int_{t}^{+\infty} \widetilde{\varphi}(x) dx$. Lemma \ref{lem:sharp4} ensures that $\mathcal{H}^\prime(t)\leqslant0$ for all $t> \max\{t_0, \beta(\alpha-\lambda)\}$. As a consequence, for all $t_1> \max\{t_0,\beta(\alpha-\lambda)\}$ and $t\geqslant t_1$,
$\mathcal{H}(t)\leqslant\mathcal{H}(t_1),$
and thus
\begin{equation}
\mathcal{E}(t) \leqslant \mathcal{E}\left(t_{1}\right)\left(\frac{t_{1}+\beta(\lambda-\alpha)}{t+\beta(\lambda-\alpha)}\right)^{\lambda \gamma-2} e^{-\widetilde{\Phi}\left(t_{1}+\beta(\lambda-\alpha)\right)+\widetilde{\Phi}(t+\beta(\lambda-\alpha))}.
\label{eq:sharp_bound}
\end{equation}
By choosing $t_1=t_0+\beta(\alpha-\lambda)$, this inequality ensures that for all $t\geqslant t_0+\beta(\alpha-\lambda)$,
\begin{equation}
\mathcal{E}(t) \leqslant \mathcal{E}\left(t_{0}+\beta(\alpha-\lambda)\right)\left(\frac{t_0}{t+\beta(\lambda-\alpha)}\right)^{\lambda \gamma-2} e^{-\widetilde{\Phi}(t_0)+\widetilde{\Phi}(t+\beta(\lambda-\alpha))}.
\end{equation}
Observe that the primitive $\widetilde{\Phi}(t)=-\int_{t}^{+\infty} \widetilde{\varphi}(x) dx$ of $\widetilde{\varphi}$ has the following expression:
\begin{equation}
    \widetilde{\Phi}(t)=-\frac{K(\alpha)}{\mu}\left(\frac{\sqrt{\mu}(1+C_0)}{t}+\frac{\alpha(1+\sqrt{2})}{(\gamma+2) t^2}+\frac{4 \alpha^{2}}{3(\gamma+2)^{2} \sqrt{\mu} t^{3}}\right),\label{eq:Phi_diff}
\end{equation}
showing that $\widetilde{\Phi}$ is non-positive. As a consequence, for all $t\geqslant t_0+\beta(\alpha-\lambda)$,
\begin{equation}
F(x(t))-F^*\leqslant e^{-\widetilde{\Phi}(t_0)}\mathcal{E}\left(t_{0}+\beta(\alpha-\lambda)\right)\frac{t_0^{\lambda\gamma-2}}{\left(t+\beta(\lambda-\alpha)\right)^{\frac{2\alpha\gamma}{\gamma+2} }},
\end{equation}
which proves the first claim of the theorem.

The value of $t_1$ can be parametrized to ensure a tight control on the energy $\mathcal{E}(t)$ in \eqref{eq:sharp_bound}. In this proof, $t_1$ is chosen as a minimizer of the following function,
$$u \mapsto (u+\beta(\lambda-\alpha))^{\lambda \gamma-2}e^{-\Phi\left(u+\beta(\lambda-\alpha)\right)}.$$
As a consequence, $u=t_1+\beta(\lambda-\alpha)$ satisfies:
\begin{equation}
\dfrac{\gamma \lambda-2}{u}  = \widetilde{\varphi}(u).
\label{eq:sharp_phi}
\end{equation}
Noticing that $\lambda \gamma-2=\frac{\gamma+2}{\alpha} K(\alpha)$, \eqref{eq:sharp_phi} can be rewritten as:
$$
\frac{\gamma+2}{\alpha u}=\frac{1}{\mu u^{2}}\left(\sqrt{\mu}(1+C_0)+\frac{2 \alpha}{(\gamma+2) u}(1+\sqrt{2})+\frac{4 \alpha^{2}}{(\gamma+2)^{2} \sqrt{\mu} u^{2}}\right)
$$
Introducing $r=(\gamma+2) \frac{\sqrt{\mu}}{\alpha} u$, this is equivalent to:
$$
r^{3}-(1+C_0)r^{2}-2(1+\sqrt{2}) r-4=0.
$$
For any $C_0>0$, the polynomial $r \mapsto r^{3}-(1+C_0)r^{2}-2(1+\sqrt{2}) r-4$ has a unique real positive root denoted $r^*$.
Defining $t_1=\frac{\alpha}{(\gamma+2) \sqrt{\mu}} r^{*}+\beta(\alpha-\lambda)$, if $t_1> \max\{t_0,\beta(\alpha-\lambda)\}$ which is guaranteed if $t_1\geqslant t_0+\beta(\alpha-\lambda)$, then the control on the energy is given by:
$$
\forall t \geqslant t_{1}, \mathcal{E}(t) \leqslant \mathcal{E}\left(t_{1}\right)\left(\frac{\alpha r^*}{\left(t+\beta(\lambda-\alpha)\right)(\gamma+2)\sqrt{\mu}}\right)^{\lambda \gamma-2} e^{-\widetilde{\Phi}\left(\frac{\alpha}{(\gamma+2) \sqrt{\mu}} r^{*}\right)+\widetilde{\Phi}(t+\beta(\lambda-\alpha))}
$$
Let $E_m$ be an energy function defined for all $t\geqslant t_0$ by:
$$E_{m}(t)=\left(1+\dfrac{\beta\alpha}{t}\right)\left(F(x(t))-F^{*}\right)+\dfrac{1}{2}\left\|\dot{x}(t)+\beta \nabla F(x(t))\right\|^{2}$$
Note that this energy is non-increasing since: $$E^{\prime}_{m}(t)=-\dfrac{\beta\alpha}{t^2}\left(F(x(t))-F^{*}\right)-\frac{\alpha}{t}\|\dot{x}(t)\|^{2}-\beta\|\nabla F(x(t))\|^{2} \leqslant 0$$
Hence, $E_{m}$ is uniformly bounded on $\left[t_{0},+\infty[\right.$. We then have:
$$
\begin{aligned}
\mathcal{E}\left(t_{1}\right) &=\left(t_{1}^{2}+t_{1} \beta(\lambda-\alpha)\right)\left(F\left(x\left(t_{1}\right)\right)-F^*\right)\\
&+\frac{1}{2}\left\|\lambda\left(x\left(t_{1}\right)-x^{*}\right)+t_{1}\left(\dot{x}\left(t_{1}\right)+\beta \nabla F\left(x\left(t_{1}\right)\right)\right)\right\|^{2} \\
&=\left(t_{1}^{2}+t_{1} \beta(\lambda-\alpha)\right) \left(F\left(x\left(t_{1}\right)\right)-F^*\right)+\frac{t_1^2}{2}\left\|\dot{x}\left(t_{1}\right)+\beta \nabla F\left(x\left(t_{1}\right)\right)\right\|^{2}\\
&+\frac{\lambda^2}{2}\left\|x\left(t_{1}\right)-x^{*}\right\|^{2}+\lambda t_{1}\left\langle\left(x\left(t_{1}\right)-x^{*}\right), \dot{x}\left(t_{1}\right)+\beta \nabla F\left(x\left(t_{1}\right)\right)\right\rangle \\
& \leqslant \left(t_{1}^{2}+t_{1} \beta(\lambda-\alpha)\right)\left(F(x(t))-F^*\right)+\frac{\lambda}{2}\left(\lambda+t_1\sqrt{\mu}\right)\|x(t_1)-x^*\|^2 \\&+\frac{t_1}{2}\left(t_1+\frac{\lambda}{\sqrt{\mu}}\right)\|\dot{x}(t)+\beta \nabla F(x(t))\|^{2}, 
\end{aligned}$$
using the inequality \begin{equation}\left|\langle x(t)-x^{*}, \dot{x}(t)+\beta \nabla F(x(t))\rangle\right| \leqslant \frac{\sqrt{\mu}}{2}\left\|x(t)-x^{*}\right\|^{2}+\frac{1}{2 \sqrt{\mu}}\|(\dot{x}(t)+\beta \nabla F(x(t)))\|^{2}.\label{eq:sharp_inmu}\end{equation} As $F$ satisfies the assumption $\mathcal{G}_\mu^2$ and noticing that $\frac{\lambda}{\sqrt{\mu}}=\frac{2}{r^*}(t_1+\beta(\lambda-\alpha))$ we get that:
$$
\begin{aligned}
\mathcal{E}(t_1)& \leqslant \left(t_{1}^{2}+t_{1} \beta(\lambda-\alpha)+\frac{\lambda^2}{\mu}+t_1\frac{\lambda}{\sqrt{\mu}}\right)\left(F(x(t))-F^*\right)\\&+\frac{t_1}{2}\left(t_1+\frac{\lambda}{\sqrt{\mu}}\right)\|\dot{x}(t)+\beta \nabla F(x(t))\|^{2} \\
& \leqslant \left(1+\frac{2}{r^*}\right)^2\left(t_{1}^{2}+t_{1} \beta(\lambda-\alpha)\right)\left(F(x(t))-F^*\right)\\&+\frac{1}{2}\left(\left(1+\frac{2}{r^*}\right)t_1^2+\frac{2}{r^*}t_1(t_1+\beta(\lambda-\alpha))\right)\|\dot{x}(t)+\beta \nabla F(x(t))\|^{2} \\
&\leqslant \left(1+\frac{2}{r^*}\right)^2t_1^2\left( F(x(t))-F^*+\frac{1}{2}\|\dot{x}(t)+\beta \nabla F(x(t))\|^{2}\right)\\
&\leqslant \left(1+\frac{2}{r^*}\right)^2t_1^2 E_m(t_1)\leqslant \left(1+\frac{2}{r^*}\right)^2t_1^2 E_m(t_0).
\end{aligned}
$$

Note that $\widetilde{\Phi}$ given in \eqref{eq:Phi_diff} is non-positive for all $t\geqslant0$ and as $t_1=\frac{\alpha}{(\gamma+2) \sqrt{\mu}} r^{*}+\beta(\alpha-\lambda)$,
$$\widetilde{\Phi}(t_1+\beta(\lambda-\alpha))=-\dfrac{\gamma+2}{\alpha} K(\alpha) \left(\frac{1+C_0}{r^{*}}+\frac{1+\sqrt{2}}{r^{* 2}}+\frac{4}{3 r^{* 3}}\right).$$

Therefore, for all $t\geqslant t_1$:
$$
F(x(t))-F^* \leqslant C_{1} e^{\frac{2 \gamma}{\gamma+2} C_{2}\left(\alpha-1-\frac{2}{\gamma}\right)}\left(1+\tfrac{\beta(\alpha-\lambda)(\gamma+2)\sqrt{\mu}}{\alpha r^*}\right) E_{m}\left(t_{0}\right)\left(\tfrac{\alpha r^*}{(\gamma+2)\sqrt{\mu}\left(t+\beta(\lambda-\alpha)\right)}\right)^{\frac{2 \alpha \gamma}{\gamma+2}},
$$
where
$$
C_{1}=\left(1+\frac{2}{r^{*}}\right)^{2},~ 
C_{2}=\frac{1+C_0}{r^{*}}+\frac{1+\sqrt{2}}{r^{* 2}}+\frac{4}{3 r^{* 3}}.
$$
Let $\mathcal{G}$ be defined as follows:
$$\mathcal{G}:t\mapsto \mathcal{H}(t)+\beta\int_{t_0+\beta(\alpha-\lambda)}^tu(u+\beta(\lambda-\alpha))^{\gamma\lambda-1}e^{-\widetilde{\Phi}(u+\beta(\lambda-\alpha))}\|\nabla F(x(u))\|^2 du,$$
Lemma \ref{lem:sharp4} guarantees that $\mathcal{G}^\prime(t)\leqslant0$ for all $t\geqslant t_0+\beta(\alpha-\lambda)$. As a consequence, for all $t\geqslant t_0+\beta(\alpha-\lambda)$,
$$\mathcal{G}(t) \leqslant \mathcal{G}(t_0+\beta(\alpha-\lambda)),$$
and as $\mathcal{H}$ is positive:
$$\beta\int_{t_0+\beta(\alpha-\lambda)}^tu(u+\beta(\lambda-\alpha))^{\gamma\lambda-1}e^{-\widetilde{\Phi}(u+\beta(\lambda-\alpha))}\|\nabla F(x(u))\|^2 du \leqslant \mathcal{G}(t_0+\beta(\alpha-\lambda)).$$
Moreover, $\widetilde{\Phi}$ is non-positive and thus:
$$\beta\int_{t_0+\beta(\alpha-\lambda)}^tu(u+\beta(\lambda-\alpha))^{\gamma\lambda-1}\|\nabla F(x(u))\|^2 du \leqslant \mathcal{G}(t_0+\beta(\alpha-\lambda)).$$
We can deduce that:
$$\int_{t_0+\beta(\alpha-\lambda)}^{+\infty}u(u+\beta(\lambda-\alpha))^{\gamma\lambda-1}\|\nabla F(x(u))\|^2 du <+\infty.$$
Note that as $u\mapsto\left(1+\beta\frac{\lambda-\alpha}{u}\right)^{\lambda\gamma-1}$ is decreasing on $(t_0+\beta(\alpha-\lambda),+\infty)$, we have that:
\begin{equation}
    \begin{aligned}
    \int_{t_0+\beta(\alpha-\lambda)}^{+\infty}u&(u+\beta(\lambda-\alpha))^{\gamma\lambda-1}\|\nabla F(x(u))\|^2 du\geqslant \\&\left(1+\beta\frac{\lambda-\alpha}{t_0+\beta(\alpha-\lambda)}\right)^{\lambda\gamma-1}\int_{t_0+\beta(\alpha-\lambda)}^{+\infty}u^{\frac{2\alpha\gamma}{\gamma+2}}\|\nabla F(x(u))\|^2 du.
    \end{aligned}
\end{equation}
In addition, the function defined by $u\mapsto u^{\frac{2\alpha\gamma}{\gamma+2}}\|\nabla F(x(u))\|^2$ is bounded on $(t_0,t_0+\beta(\alpha-\lambda))$ and consequently:
\begin{equation}
    \int_{t_0}^{+\infty}u^{\frac{2\alpha\gamma}{\gamma+2}}\|\nabla F(x(u))\|^2 du <+\infty.
\end{equation}
\textbf{Case $\alpha=1+\frac{2}{\gamma}$ (Proof of statements 2 and 3).}\\
Lemma \ref{lem:sharp3} ensures that for all $t> \max\{t_0,\beta\}$,
\begin{equation}
\mathcal{E}^\prime(t)+\beta t(t-\beta)\|\nabla F(x(t))\|^2\leqslant\frac{\beta}{t(t-\beta)}\mathcal{E}(t),
\label{eq:sharp_alpha=}
\end{equation}
noticing that $\alpha-\lambda=1$. This inequality implies that $t\mapsto \mathcal{E}(t)e^{\frac{\beta}{t-\beta}}$ is decreasing on $(t_0+\beta,+\infty)$. Consequently, for all $t\geqslant t_0+\beta$,
\begin{equation*}
\mathcal{E}(t)\leqslant\mathcal{E}(t_0+\beta)e^{-\frac{\beta}{t-\beta}+\frac{\beta}{t_0}}\leqslant \mathcal{E}(t_0+\beta)e^{\frac{\beta}{t_0}}.
\end{equation*}
Moreover, 
\begin{equation*}
\begin{aligned}
\mathcal{E}(t_0+\beta)&=t_0\left(t_0+\beta\right)\left(F(x(t_0+\beta))-F^*\right)\\&+\frac{1}{2}\|\lambda(x(t_0+\beta)-x^*)+t(\dot x(t_0+\beta)+\beta\nabla F(x(t_0+\beta)))\|^2\\
&\leqslant \left((t_0+\beta)^2+\frac{\lambda^2+\sqrt{\mu}}{\mu}\right)\left(F(x(t_0+\beta))-F^*\right)\\&+\frac{(t_0+\beta)^2+\frac{1}{\sqrt{\mu}}}{2}\|\dot x(t_0+\beta)+\beta\nabla F(x(t_0+\beta))\|^2\\
&\leqslant \left((t_0+\beta)^2+\frac{\lambda^2+\sqrt{\mu}}{\mu}\right)E_m(t_0+\beta)
\leqslant\left((t_0+\beta)^2+\frac{\lambda^2+\sqrt{\mu}}{\mu}\right)E_m(t_0),
\end{aligned}
\end{equation*}
using inequality \eqref{eq:sharp_inmu}. Hence, for all $t\geqslant t_0+\beta$,
\begin{equation}
F(x(t))-F^*\leqslant \left((t_0+\beta)^2+\frac{\lambda^2+\sqrt{\mu}}{\mu}\right)e^{\frac{\beta}{t_0}}\frac{E_m(t_0)}{t(t-\beta)}.
\end{equation}
Inequality \eqref{eq:sharp_alpha=} also guarantees that $$t\mapsto \mathcal{E}(t)e^{\frac{\beta}{t-\beta}}+\int_{t_0+\beta}^t\beta u(u-\beta)e^{\frac{\beta}{u-\beta}}\|\nabla F(x(u))\|^2du,$$ is bounded on $(t_0+\beta,+\infty)$. As $\mathcal{E}(t)e^{\frac{\beta}{t-\beta}}$ is positive for all $t\geqslant t_0+\beta$, we can deduce that there exists $M>0$ such that for all $t\geqslant t_0+\beta$,
\begin{equation*}
\int_{t_0+\beta}^t (u-\beta)^2\|\nabla F(x(u))\|^2du\leqslant \int_{t_0+\beta}^t u(u-\beta)e^{\frac{\beta}{u-\beta}}\|\nabla F(x(u))\|^2du<M,
\end{equation*}
and thus,
\begin{equation}
\int_{t_0+\beta}^{+\infty} (u-\beta)^2\|\nabla F(x(u))\|^2du<+\infty.
\end{equation}
By using the same arguments as in the first case, we can conclude that:
\begin{equation}
\int_{t_0}^{+\infty} u^2\|\nabla F(x(u))\|^2du<+\infty.
\end{equation}

\qed
\subsection{Proof of Theorem \ref{thm:sharp2}}

\label{sec:proof_sharp2}

Let $F$ be a convex $C^2$ function satisfying $\mathcal{G}_\mu^2$ for some $\mu>0$. The convexity of $F$ implies that $F$ satisfies $\mathcal{H}_1$ and the following lemma ensures that $F$ also satisfies $\mathcal{H}_{2-\delta}^{loc}$ for all $\delta\in(0,1]$. The proof of this lemma is given in Section \ref{sec:proof_H1loc}.

\begin{lemme}
Let $F: \mathbb{R}^{n} \rightarrow \mathbb{R}$ be a convex $C^2$ function with a non empty set of minimizer $X^*$. Then, for all $\delta\in(0,1]$, the function $F$ satisfies $\mathcal{H}_{2-\delta}^{loc}$.
\label{lem:H1loc}
\end{lemme}
Let $\alpha\geqslant3$, $\beta>0$ and $\varepsilon\in(0,1)$. As $F$ satisfies $\mathcal{H}_{1}$, the first and second claims of Theorem \ref{thm:sharp1} ensure that there exists a decreasing function $\phi$ such that:
\begin{equation*}
\forall t\geqslant t_0+\frac{\beta\alpha}{3},\quad F(x(t))-F^*\leqslant \phi(t),
\end{equation*}
where $\phi(t)\rightarrow0$ as $t\rightarrow+\infty$. Therefore, as $F$ satisfies $\mathcal{G}_\mu^2$, for all $\nu>0$, there exists $T\geqslant t_0+\frac{\beta\alpha}{3}$ such that for all $t\geqslant T$, $x(t)\in B(x^*,\nu)$. As a consequence, for all $\delta\in(0,1]$, there exists $T\geqslant t_0+\frac{\beta\alpha}{3}$ such that for all $t\geqslant T$,
\begin{equation*}
F(x(t))-F^*\leqslant \frac{1}{2-\delta}\left\langle x(t)-x^*,\nabla F(x(t))\right\rangle.
\end{equation*}
Let $\delta=\frac{4\varepsilon}{\alpha+\varepsilon}$. As $\alpha\geqslant3$ and $\varepsilon\in(0,1)$, the condition $\alpha\geqslant 1+\frac{2}{2-\delta}$ is satisfied. Then, by setting $T$ as the initial time in \eqref{eq:Hessian_ODE}, the first claim of Theorem \ref{thm:sharp1} gives the first result and the third claim of Theorem \ref{thm:sharp1} guarantees that there exists $M>0$ such that
\begin{equation*}
\int_{T}^{+\infty}u^{\alpha-\varepsilon}\|\nabla F(x(u))\|^2du<M.
\end{equation*}
As $F$ satisfies $\mathcal{G}_\mu^2$, Lemma \ref{lem:Loja} ensures that 
\begin{equation*}
\int_{T}^{+\infty}u^{\alpha-\varepsilon}\left(F(x(u))-F^*\right) du <\frac{M}{2\mu}.
\end{equation*}
On the other hand, as $u\mapsto u^{\alpha-\varepsilon}(F(x(u))-F^*)$ is bounded on $\left(t_0,T\right)$, we have that:
\begin{equation*}
\int_{t_0}^{T}u^{\alpha-\varepsilon}\left(F(x(u))-F^*\right) du <+\infty,
\end{equation*}
and consequently,
\begin{equation*}
\int_{t_0}^{+\infty}u^{\alpha-\varepsilon}\left(F(x(u))-F^*\right) du <+\infty.
\end{equation*}
\qed

\subsection{Proof of Theorem \ref{thm:flat}}
\label{sec:proof_flat1}

We define $\mathcal{E}$ as the following Lyapunov function:
$$
\begin{aligned}
\mathcal{E}(t) =&\left(t^{2}+t \beta(\lambda-\alpha)\right)\left(F(x(t))-F^{*}\right)+\frac{\xi}{2}\left\|x(t)-x^{*}\right\|^{2}\\
&+\frac{1}{2}\left\|\lambda\left(x(t)-x^{*}\right)+t(\dot{x}(t)+\beta \nabla F(x(t)))\right\|^{2},
\end{aligned}
$$
where $x^*$ is the unique minimizer of $F$, $\lambda\in\R$ and $\xi\in\R$. Let $\mathcal{H}$ be the function defined as follows
$$
\mathcal{H}(t)=t^{p} \mathcal{E}(t),
$$
where $p>0$. Using the notations
$$
\begin{aligned}
&a(t)=t\left(F(x(t))-F^{*}\right), \quad b(t)=\frac{1}{2t}\left\|\lambda\left(x(t)-x^{*}\right)+t( \dot{x}(t)+\beta \nabla F(x(t)))\right\|^{2}, \\
&c(t)=\frac{1}{2t}\left\|x(t)-x^{*}\right\|^{2},
\end{aligned}
$$
we have
$$
\mathcal{E}(t)=(t+\beta(\lambda-\alpha))a(t)+t(b(t)+\xi c(t)).
$$
Let $p=\frac{4}{\gamma_{1}-2}$,  $\lambda=\frac{2}{\gamma_{1}-2}$ and $\xi=\lambda(\lambda+1-\alpha)$.

\begin{lemme}
Let $p=\frac{4}{\gamma_{1}-2}$,  $\lambda=\frac{2}{\gamma_{1}-2}$ and $\xi=\lambda(\lambda+1-\alpha)$, for all $t \geqslant \max(t_0,\beta (\alpha-\lambda),\beta (2(\alpha-\lambda)-1)))$
\begin{equation}
\begin{aligned}\mathcal{E}^{\prime}(t) &\leqslant((2-\lambda\gamma_1)t+\beta(\lambda-\alpha-\lambda\gamma_1(2(\lambda-\alpha)+1))))\left(F(x(t))-F^{*}\right) \\ 
& +2(\lambda+1-\alpha)b(t)-2\lambda^2(\lambda+1-\alpha)c(t)-\beta t(t+\beta (\lambda-\alpha)) \| \nabla F(x(t))\|^{2}. \\
\end{aligned}
\end{equation}\label{lem:flat2}
The proof of this lemma is given in Section \ref{sec:proof_flat2}
\end{lemme}
\noindent Consequently, for all $t \geqslant \max(t_0,\beta (\alpha-\lambda),\beta (2(\alpha-\lambda)-1)))$:
$$
\begin{aligned}
\mathcal{H}^{\prime}(t)&=t^{p-1}\left(p \mathcal{E}(t)+t \mathcal{E}^{\prime}(t)\right)\leqslant t^{p-1}[p \mathcal{E}(t)+2t(\lambda+1-\alpha)b(t)\\
&+t((2-\lambda\gamma_1)t+\beta(\lambda-\alpha-\lambda\gamma_1(2(\lambda-\alpha)+1))\left(F(x(t))-F^{*}\right)\\
&-2t(\lambda^2(\lambda+1-\alpha))c(t)-\beta t^2(t+\beta(\lambda-\alpha))\|\nabla F(x(t))\|^{2}]\\ 
&\leqslant t^{p}((2-\gamma_1 \lambda+p) a(t)+2 (\lambda-\alpha+1)+p) b(t)+\lambda(\lambda+1-\alpha)(p-2 \lambda) c(t))\\
&+t^{p-1}\beta((p+1)(\lambda-\alpha)-\lambda\gamma_1(2(\lambda-\alpha)+1))a(t)\\&-\beta t^{p+1}(t+\beta(\lambda-\alpha))\|\nabla F(x(t))\|^{2}.
\end{aligned}
$$
As $p=\frac{4}{\gamma_{1}-2}$ and $\lambda=\frac{2}{\gamma_{1}-2}$ this implies that
\begin{equation}
\mathcal{H}^{\prime}(t)\leqslant 2 t^{p}\left(\frac{\gamma_{1}+2}{\gamma_{1}-2}-\alpha\right)b(t)+t^{p-1}\beta C_1 a(t)-\beta t^{p+1}(t+\beta(\lambda-\alpha))\|\nabla F(x(t))\|^{2},
\label{eq:Hprime1}
\end{equation}
where $C_1=(p+1)(\lambda-\alpha)-\lambda\gamma_1(2(\lambda-\alpha)+1)$. Under the assumption $\alpha\geqslant\frac{\gamma_1+2}{\gamma_1-2}$, $C_1$ is strictly positive and \eqref{eq:Hprime1} ensures that\begin{equation}
\mathcal{H}^{\prime}(t)\leqslant t^{p-1}\beta C_1 a(t)-\beta t^{p+1}(t+\beta(\lambda-\alpha))\|\nabla F(x(t))\|^{2}.
\label{eq:Hprime2}
\end{equation}
We define $\mathcal{G}$ as follows:
$$\mathcal{G}:t\mapsto \mathcal{H}(t)-\beta C_1 \int_{t_1}^t u^{p-1}a(u)du,$$
where $t_1>\max\left\{t_0,\beta(2(\alpha-\lambda)-1),t_m \right\}$ and $t_m>\beta(\alpha-\lambda)$ satisfies
\begin{equation}\frac{t_m^p}{(t_m+\beta(\lambda-\alpha))^{p+1}}\beta C_1\leqslant\frac{1}{2}.\end{equation}
As $t\mapsto\frac{t^p}{(t+\beta(\lambda-\alpha))^{p+1}}$ is decreasing on $\left(\beta(\alpha-\lambda),+\infty\right)$ and tends towards $0$, $t_m$ is well defined. Equation \eqref{eq:Hprime2} implies that $\mathcal{G}^\prime(t)\leqslant0$ for all $t\geqslant t_1$ and therefore there exists $A\in\R$ such that $\mathcal{G}(t)\leqslant A$ for all $t\geqslant t_1$, and for all $t\geqslant t_1$,
$$\begin{aligned}
\mathcal{G}(t)&=t^p((t+\beta(\lambda-\alpha))a(t)+tb(t)+t\xi c(t))-\beta C_1 \int_{t_1}^t u^{p-1}a(u)du\\
&\geqslant t^p((t+\beta(\lambda-\alpha))a(t)+t\xi c(t))-\beta C_1 \int_{t_1}^t u^{p-1}a(u)du\\
&\geqslant (t+\beta(\lambda-\alpha))^{p+1}a(t)+t^{p+1}\xi c(t)-\beta C_1 \int_{t_1}^t u^{p-1}a(u)du.\\
\end{aligned}$$
Moreover, $-\beta C_1 \int_{t_1}^t u^{p-1}a(u)du\geqslant-\left(\frac{t_1}{t_1+\beta(\lambda-\alpha)}\right)^{p+1}\beta C_1\int_{t_1}^t\frac{(u+\beta(\lambda-\alpha))^{p+1}a(u)}{u^2}du$. Recall that $F$ satisfies $\mathcal{G}_\mu^{\gamma_2}$, thus there exists $K=\frac{\mu}{2}>0$ such that
\begin{equation*}
\forall x\in\R^n,\quad K d\left(x, X^{*}\right)^{\gamma_{2}} \leqslant F(x)-F^*,
\end{equation*}
and therefore
$$\begin{aligned}
t^{p+1}\xi c(t)&=t^p\frac{\xi}{2}\|x(t)-x^*\|^2=t^p\frac{\xi}{2K^{\frac{2}{\gamma_2}}}\left(K\|x(t)-x^*\|^{\gamma_2}\right)^{\frac{2}{\gamma_2}}.
\end{aligned}$$
As $\xi<0$ and $F$ has a unique minimizer,
$$\begin{aligned}
t^{p+1}\xi c(t)&=t^p\frac{\xi}{2K^{\frac{2}{\gamma_2}}}\left(Kd(x(t),X^*)^{\gamma_2}\right)^{\frac{2}{\gamma_2}} \geqslant t^p\frac{\xi}{2K^{\frac{2}{\gamma_2}}}\left(F(x(t))-F^*\right)^{\frac{2}{\gamma_2}}\\
&\geqslant t^{p-\frac{2}{\gamma_2}}\frac{\xi}{2K^{\frac{2}{\gamma_2}}}a(t)^{\frac{2}{\gamma_2}}\geqslant t^{p-\frac{2}{\gamma_2}-(p+1)\frac{2}{\gamma_2}}\frac{\xi}{2K^{\frac{2}{\gamma_2}}}\left(t^{p+1}a(t)\right)^{\frac{2}{\gamma_2}}\\
&\geqslant t^{p-\frac{2}{\gamma_2}-(p+1)\frac{2}{\gamma_2}}\left(\frac{t_1}{t_1+\beta(\lambda-\alpha)}\right)^{\frac{2(p+1)}{\gamma_2}}\frac{\xi}{2K^{\frac{2}{\gamma_2}}}\left((t+\beta(\alpha-\lambda))^{p+1}a(t)\right)^{\frac{2}{\gamma_2}}.
\end{aligned}$$
Recall that $p=\frac{4}{\gamma_1-2}$, therefore $p-\frac{2}{\gamma_2}-(p+1)\frac{2}{\gamma_2}=\frac{4(\gamma_2-\gamma_1)}{\gamma_2(\gamma_1-2)}\leqslant 0$ since $\gamma_1\geqslant\gamma_2$. As $t\geqslant t_1$, we have
$$
t^{p+1}\xi c(t)\geqslant \frac{t_1^{p-\frac{2}{\gamma_2}}}{(t_1+\beta(\lambda-\alpha))^{\frac{2(p+1)}{\gamma_2}}}\frac{\xi}{2K^{\frac{2}{\gamma_2}}}\left((t+\beta(\alpha-\lambda))^{p+1}a(t)\right)^{\frac{2}{\gamma_2}}.
$$
We define $v:t\mapsto(t+\beta(\lambda-\alpha))^{p+1}a(t)$ for all $t\geqslant t_1$. Then, for all $t\geqslant t_1$:
\begin{equation}
\mathcal{G}(t)\geqslant v(t)-C_2v(t)^{\frac{2}{\gamma_2}}-\left(\frac{t_1}{t_1+\beta(\lambda-\alpha)}\right)^{p+1}\beta C_1\int_{t_1}^t\frac{v(u)}{u^2}du,
\label{eq:control_Gv}
\end{equation}
where $C_2=-\frac{t_1^{p-\frac{2}{\gamma_2}}}{(t_1+\beta(\lambda-\alpha))^{\frac{2(p+1)}{\gamma_2}}}\frac{\xi}{2K^{\frac{2}{\gamma_2}}}>0$.
Let $t_2>t_1$ and $t^*=\underset{t\in[t_1,t_2]}{\text{argmax}}~v(t)$. Then, 
$$\begin{aligned}
\mathcal{G}(t^*)&\geqslant v(t^*)-C_2v(t^*)^{\frac{2}{\gamma_2}}-\left(\frac{t_1}{t_1+\beta(\lambda-\alpha)}\right)^{p+1}\beta C_1\int_{t_1}^{t^*}\frac{v(u)}{u^2}du\\
&\geqslant v(t^*)-C_2 v(t^*)^{\frac{2}{\gamma_2}}-\left(\frac{t_m}{t_m+\beta(\lambda-\alpha)}\right)^{p+1}\beta C_1\int_{t_m}^{+\infty}\frac{v(t^*)}{u^2}du\\
&\geqslant v(t^*)-C_2 v(t^*)^{\frac{2}{\gamma_2}}-\frac{t_m^p}{\left(t_m+\beta(\lambda-\alpha)\right)^{p+1}}\beta C_1v(t^*)\\
&\geqslant \frac{1}{2}v(t^*)-C_2 v(t^*)^{\frac{2}{\gamma_2}}.
\end{aligned}$$
As $\mathcal{G}(t)\leqslant A$ for all $t\geqslant t_1$, we get that:
$$v(t^*)-2C_2 v(t^*)^{\frac{2}{\gamma_2}}\leqslant 2A,$$
and consequently
\begin{equation}
v(t^*)^{\frac{2}{\gamma_2}}\left(v(t^*)^{1-\frac{2}{\gamma_2}}-2C_2\right)\leqslant 2A.
\label{eq:control_v1}
\end{equation}
\begin{lemme}
Let $x\in\R^+$, $\delta\in(0,1)$, $K_1>0$ and $K_2>0$. Then,
$$x^\delta(x^{1-\delta}-K_1)\leqslant K_2\quad \implies\quad x\leqslant \left(K_2^{1-\delta}+K_1\right)^{\frac{1}{1-\delta}}.$$
\label{lem:control_v}
\end{lemme}
Applying Lemma \ref{lem:control_v} to \eqref{eq:control_v1} we get that 
\begin{equation}
v(t^*)\leqslant\left((2A)^{1-\frac{2}{\gamma_2}}+2C_2\right)^{\frac{\gamma_2}{\gamma_2-2}},
\end{equation}
and thus for all $t\in[t_1,t_2]$
\begin{equation}
v(t)\leqslant\left((2A)^{1-\frac{2}{\gamma_2}}+2C_2\right)^{\frac{\gamma_2}{\gamma_2-2}}.
\end{equation}
This bound does not depend on $t_2$ so we can deduce that $v$ is bounded on $[t_1,+\infty)$.
As a consequence, there exists $M>0$ such that for all $t\geqslant t_1$:
$$v(t)=(t+\beta(\lambda-\alpha))^{p+1}a(t)=t(t+\beta(\lambda-\alpha))^{p+1}(F(x(t))-F^*)\leqslant M,$$
which implies that
\begin{equation}
F(x(t))-F^*\leqslant \frac{M}{t(t+\beta(\lambda-\alpha))^{p+1}}\leqslant \left(\frac{t_1}{t_1+\beta(\lambda-\alpha)}\right)^{p+1} \frac{M}{t^{p+2}},
\end{equation}
i.e. as $t\rightarrow+\infty$
\begin{equation}
F(x(t))-F^*=\mathcal{O}\left(t^{-\frac{2\gamma_1}{\gamma_1-2}}\right).
\end{equation}
Let $\mathcal{F}$ be defined by
$$\mathcal{F}:t\mapsto \mathcal{H}(t)-\beta C_1 \int_{t_1}^t u^{p-1}a(u)du+\beta \int_{t_1}^t u^{p+1}(u+\beta(\lambda-\alpha))\|\nabla F(x(u))\|^2du.$$
Equation \eqref{eq:Hprime2} implies that $\mathcal{F}^\prime(t)\leqslant0$ for all $t\geqslant t_1$ and therefore there exists $B\in\R$ such that $\mathcal{F}(t)\leqslant B$ for all $t\geqslant t_1$. By applying \eqref{eq:control_Gv} we get that for all $t\geqslant t_1$
$$\begin{aligned}\mathcal{F}(t)&\geqslant v(t)-C_2v(t)^{\frac{2}{\gamma_2}}-\left(\frac{t_1}{t_1+\beta(\lambda-\alpha)}\right)^{p+1}\beta C_1\int_{t_1}^t\frac{v(u)}{u^2}du\\&+\beta \int_{t_1}^t u^{p+1}(u+\beta(\lambda-\alpha))\|\nabla F(x(u))\|^2du.\end{aligned}$$
We proved that there exists $M>0$ such that for all $t\geqslant t_1$, $v(t)\leqslant M$. Hence,
$$-\left(\frac{t_1}{t_1+\beta(\lambda-\alpha)}\right)^{p+1}\beta C_1\int_{t_1}^t\frac{v(u)}{u^2}du\geqslant -M\beta C_1\frac{t_1^p}{(t_1+\beta(\lambda-\alpha))^{p+1}}.$$
\begin{lemme}
Let $g:x\mapsto x-Kx^\delta$ for some $K>0$ and $\delta\in(0,1)$. Then for all $x\geqslant 0$,
$$g(x)\geqslant K(\delta-1)(\delta K)^{\frac{\delta}{1-\delta}}.$$
\label{lem:min_g}
\end{lemme}
Lemma \ref{lem:min_g} ensures that for all $t\geqslant t_1$
\begin{equation}
v(t)-C_2v(t)^{\frac{2}{\gamma_2}}\geqslant -C_2\left(1-\frac{2}{\gamma_2}\right)\left(\frac{2C_2}{\gamma_2}\right)^\frac{2}{\gamma_2-2}.
\end{equation}
Thus, 
$$\begin{aligned}
\mathcal{F}(t)&\geqslant -C_2\left(1-\frac{2}{\gamma_2}\right)\left(\frac{2C_2}{\gamma_2}\right)^\frac{2}{\gamma_2-2}-M\beta C_1\frac{t_1^p}{(t_1+\beta(\lambda-\alpha))^{p+1}}\\&+\beta \int_{t_1}^t u^{p+1}(u+\beta(\lambda-\alpha))\|\nabla F(x(u))\|^2du.
\end{aligned}$$
As there exists $B\in\R$ such that $\mathcal{F}(t)\leqslant B$, for all $t\geqslant t_1$, we can deduce that
$$\begin{aligned}
\beta\int_{t_1}^t u^{p+1}(u+\beta(\lambda-\alpha))\|\nabla F(x(u))\|^2du&\leqslant B+C_2\left(1-\frac{2}{\gamma_2}\right)\left(\frac{2C_2}{\gamma_2}\right)^\frac{2}{\gamma_2-2}\\&+M\beta C_1\frac{t_1^p}{(t_1+\beta(\lambda-\alpha))^{p+1}},
\end{aligned}$$
and therefore
\begin{equation}
\int_{t_1}^{+\infty} (u+\beta(\lambda-\alpha))^{\frac{2\gamma_1}{\gamma_1-2}}\|\nabla F(x(u))\|^2du<+\infty.
\end{equation}
By using the same arguments as in the proof of Theorem 1 and the boundedness of $u\mapsto(u+\beta(\lambda-\alpha))^{\frac{2\gamma_1}{\gamma_1-2}}\|\nabla F(x(u))\|^2$ on $(t_0,t_1)$, we conclude that:
\begin{equation}
\int_{t_0}^{+\infty} u^{\frac{2\gamma_1}{\gamma_1-2}}\|\nabla F(x(u))\|^2du<+\infty.
\end{equation}
\qed

\appendix

\section{Appendix}

\subsection{Proof of Corollary \ref{cor:sharp2}}
\label{sec:proof_sharp3}
The first claim is obtained by applying the following lemma to Theorem \ref{thm:sharp2}. The proof of this lemma is given in Section \ref{sec:proof_delta+1}.

\begin{lemme}
Let $F: \mathbb{R}^{n} \rightarrow \mathbb{R}$ be a convex function having a non empty set of minimizers where $F^*=\inf\limits_{x\in\R^n}F(x)$. Assume that for some $t_1>0$ and $\delta>0$, $F$ satisfies:
$$\int_{t_1}^{+\infty}u^\delta(F(x(u))-F^*)du<+\infty.$$
Let $z:t\mapsto\frac{\int_{t/2}^tu^\delta x(u)du}{\int_{t/2}^tu^\delta du}$. Then, as $t\rightarrow+\infty$,
\begin{equation}
F(z(t))-F^*=o\left(t^{-\delta-1}\right).
\end{equation}
\label{lem:delta+1}
\end{lemme}
The second and third claim are proved by applying Lemma \ref{lem:delta+1phi} to $\phi:x\mapsto F(x)-F^*$. The proof of this lemma is given in Section \ref{subsection:delta+1phi}.

\begin{lemme}
Let $\phi: \mathbb{R}^{n} \rightarrow \mathbb{R}^+$ such that for some $t_1>0$ and $\delta>0$, $\phi$ satisfies:
\begin{equation*}
\int_{t_1}^{+\infty}u^\delta\phi(x(u))du<+\infty.
\end{equation*}
Then, as $t\rightarrow+\infty$,
\begin{equation}
\inf\limits_{u\in[t/2,t]}\phi(x(u))=o\left(t^{-\delta-1}\right)\quad\mbox{ and }\quad
\liminf\limits_{t\rightarrow+\infty} t^{\delta+1}\log (t)\phi(x(t))=0.
\end{equation}
\label{lem:delta+1phi}
\end{lemme}

\subsection{Proof of Corollary \ref{cor:flat1}}
\label{sec:proof_flat_cor}

Let $F: \mathbb{R}^{n} \rightarrow \mathbb{R}$ be a convex $C^2$ function having a unique minimizer $x^*$. Assume that $F$ satisfies $\mathcal{H}_{\gamma_1}$ and $\mathcal{G}_\mu^{\gamma_2}$ for some $\gamma_1>2$, $\gamma_2>2$ such that $\gamma_{1}\geqslant\gamma_{2}$ and $\mu>0$. Let $x$ be a solution of \eqref{eq:Hessian_ODE} for all $t\geqslant t_0$ where $t_0>0$, $\alpha\geqslant\frac{\gamma_{1}+2}{\gamma_{1}-2}$ and $\beta>0$. Theorem \ref{thm:flat} ensures that:
\begin{equation*}
\int_{t_0}^{+\infty}u^{\frac{2\gamma_1}{\gamma_1-2}}\|\nabla F(x(u))\|^2du<+\infty.
\end{equation*}
Moreover, as $F$ satisfies $\mathcal{G}_\mu^{\gamma_2}$ for some $\gamma_2>2$, Lemma \ref{lem:Loja} implies that:
\begin{equation}
\int_{t_0}^{+\infty}u^{\frac{2\gamma_1}{\gamma_1-2}}\left(F(x(u))-F^*\right)^{\frac{2(\gamma_2-1)}{\gamma_2}}du<+\infty.
\label{eq:cor2_1}
\end{equation}
By applying Lemma \ref{lem:delta+1phi} to $\phi:x\mapsto\left(F(x)-F^*\right)^{\frac{2(\gamma_2-1)}{\gamma_2}}$, we get that as $t$ tends to $+\infty$,
\begin{equation*}
\inf\limits_{u\in\left[t/2,t\right]}\left(F(x(u))-F^*\right)^{\frac{2(\gamma_2-1)}{\gamma_2}}=o\left(t^{-\frac{3\gamma_1-2}{\gamma_1-2}}\right).
\end{equation*}
Hence,
\begin{equation*}
\inf\limits_{u\in\left[t/2,t\right]}F(x(u))-F^*=o\left(t^{-\frac{(3\gamma_1-2)\gamma_2}{2(\gamma_1-2)(\gamma_2-1)}}\right).
\end{equation*}

\subsection{Proof of Lemma \ref{lem:sharp1}}
\label{subsection:Proof for Lemma 1}

Recall that for all $t\geqslant t_0$:
\begin{equation*}
\mathcal{E}(t) =\left(t^{2}+t \beta(\lambda-\alpha)\right)\left(F(x(t))-F^{*}\right)+\frac{1}{2}\left\|\lambda\left(x(t)-x^{*}\right)+t(\dot{x}(t)+\beta \nabla F(x(t)))\right\|^{2},
\end{equation*}
The Lyapunov function $\mathcal{E}$ is differentiable and simple calculations give that:
\begin{equation*}
\begin{aligned}
\mathcal{E}^{\prime}(t) &=(2 t+\beta(\lambda-\alpha))\left(F(x(t))-F^{*}\right)+\left(t^{2}+t \beta(\lambda-\alpha)\right)\langle\nabla F(x(t)), \dot{x}(t)\rangle \\
&+\lambda(\lambda+1-\alpha)\left\langle x(t)-x^{*}, \dot{x}(t)\right\rangle+\lambda(\beta-t)\left\langle\nabla F(x(t)), x(t)-x^{*}\right\rangle \\
&+t(\lambda+1-\alpha)\|\dot{x}(t)\|^2+t(\beta-t)\langle\nabla F(x(t)), \dot{x}(t)\rangle \\
&+t \beta(\lambda+1-\alpha)\langle\nabla F(x(t)), \dot{x}(t)\rangle+t \beta(\beta-t)\|\nabla F(x(t))\|^2 \\
&=(2 t+\beta(\lambda-\alpha))\left(F(x(t))-F^{*}\right)+2 t \beta(\lambda+1-\alpha)\langle\nabla F(x(t)), \dot{x}(t)\rangle \\
&+\lambda(\lambda+1-\alpha)\left\langle x(t)-x^{*}, \dot{x}(t)\right\rangle+\lambda(\beta-t)\left\langle\nabla F(x(t)), x(t)-x^{*}\right\rangle \\
&+t(\lambda+1-\alpha)\|\dot{x}(t)\|^2+t \beta(\beta-t)\|\nabla F(x(t))\|^2. \\
 \\
\end{aligned}
\end{equation*}
By rearranging the terms, we get that:
\begin{equation*}
\begin{aligned}
\mathcal{E}^\prime(t)&=(2 t+\beta(\lambda-\alpha))\left(F(x(t))-F^{*}\right) \\
&+\frac{\lambda+1-\alpha}{t}\left[t^{2}\|\dot{x}(t)\|^2+t^{2} \beta^{2}\|\nabla F(x(t))\|^2+2 t^{2} \beta\langle\nabla F(x(t)), \dot{x}(t)\rangle\right] \\
&+\frac{\lambda(\lambda+1-\alpha)}{t}\left[\left\langle x(t)-x^{*}, t \dot{x}(t)\right\rangle+\left\langle x(t)-x^{*}, t \beta \nabla F(x(t))\right\rangle\right] \\
&-t \beta(t+\beta(\lambda-\alpha))\|\nabla F(x(t))\|^{2}-\lambda(t+\beta(\lambda-\alpha))\left\langle\nabla F(x(t)), x(t)-x^{*}\right\rangle\\
&=(2 t+\beta(\lambda-\alpha))\left(F(x(t))-F^{*}\right)+\frac{\lambda+1-\alpha}{t}\|t(\dot{x}(t)+\beta \nabla F(x(t)))\|^{2} \\
&+\frac{\lambda(\lambda+1-\alpha)}{t}\left\langle x(t)-x^{*}, t(\dot{x}(t)+\beta \nabla F(x(t)))\right\rangle \\
&-t \beta(t+\beta(\lambda-\alpha))\|\nabla F(x(t))\|^{2}-\lambda(t+\beta(\lambda-\alpha))\left\langle\nabla F(x(t)), x(t)-x^{*}\right\rangle.\\
\end{aligned}
\end{equation*}
A last step allows us to conclude that:
\begin{equation*}
\begin{aligned}
\mathcal{E}^\prime(t)&=(2-\gamma \lambda)(t+\beta(\lambda-\alpha))\left(F(x(t))-F^{*}\right)-\beta(\lambda-\alpha)\left(F(x(t))-F^{*}\right) \\
&-\lambda(t+\beta(\lambda-\alpha))\left[-\gamma\left(F(x(t))-F^{*}\right)+\left\langle\nabla F(x(t)), x(t)-x^{*}\right\rangle\right]\\
&+\frac{\lambda(\lambda+1-\alpha)}{t}\left\langle x(t)-x^{*}, t(\dot{x}(t)+\beta \nabla F(x(t)))\right\rangle \\
&+\frac{\lambda+1-\alpha}{t}\|t(\dot{x}(t)+\beta \nabla F(x(t)))\|^{2}-t \beta(t+\beta(\lambda-\alpha))\|\nabla F(x(t))\|^{2}.
\end{aligned}
\end{equation*}
We distinguish several terms containing $F(x(t))-F^*$ as they will be treated separately.

\subsection{Proof of Lemma \ref{lem:sharp3}}
\label{subsection:Proof for Lemma 3}
Notice that for all $t\geqslant t_0$,
$$
\begin{aligned}
\frac{\gamma \lambda-2}{t} \mathcal{E}(t)=&-(2-\gamma \lambda)(t+\beta(\lambda-\alpha))\left(F(x(t))-F^{*}\right)\\
&+\frac{1}{2} \frac{\gamma \lambda-2}{t}\left\|\lambda\left(x(t)-x^{*}\right)+t(\dot{x}(t)+\beta \nabla F(x(t)))\right\|^{2}.
\end{aligned}
$$
By applying Lemma \ref{lem:sharp1}, we get that:
$$
\begin{aligned}
\mathcal{E}^{\prime}(t)+\frac{\gamma \lambda-2}{t} \mathcal{E}(t) &=\frac{1}{2} \frac{\gamma \lambda-2}{t}\left\|\lambda\left(x(t)-x^{*}\right)+t(\dot{x}(t)+\beta \nabla F(x(t)))\right\|^{2} \\
&-\lambda(t+\beta(\lambda-\alpha))\left[-\gamma\left(F(x(t))-F^{*}\right)+\left\langle\nabla F(x(t)), x-x^{*}\right\rangle\right]\\
&+\frac{\lambda(\lambda+1-\alpha)}{t}\left\langle x(t)-x^{*}, t(\dot{x}(t)+\beta \nabla F(x(t)))\right\rangle \\
&+\frac{\lambda+1-\alpha}{t}\|t(\dot{x}(t)+\beta \nabla F(x(t)))\|^{2}-t \beta(t+\beta(\lambda-\alpha))\|\nabla F(x(t))\|^{2} \\
&+\beta(\alpha-\lambda)\left(F(x(t))-F^{*}\right).
\end{aligned}
$$
As $F$ satisfies $\mathcal{H}_\gamma$, for all $t\geqslant \max\{ t_0,\beta(\alpha-\lambda)\}$
\begin{equation*}
\lambda(t+\beta(\lambda-\alpha))\left[-\gamma\left(F(x(t))-F^{*}\right)+\left\langle\nabla F(x(t)), x-x^{*}\right\rangle\right]\geqslant0,
\end{equation*}
and hence
$$
\begin{aligned}
\mathcal{E}^{\prime}(t)+\frac{\gamma \lambda-2}{t} \mathcal{E}(t) & \leqslant \frac{1}{2} \frac{\gamma \lambda-2}{t}\left\|\lambda\left(x(t)-x^{*}\right)+t(\dot{x}(t)+\beta \nabla F(x(t)))\right\|^{2} \\
&+\frac{\lambda(\lambda+1-\alpha)}{t}\left\langle x(t)-x^{*}, t(\dot{x}(t)+\beta \nabla F(x(t)))\right\rangle \\
&+\frac{\lambda+1-\alpha}{t}\|t(\dot{x}(t)+\beta \nabla F(x(t)))\|^{2} \\
&+\beta(\alpha-\lambda)\left(F(x(t))-F^{*}\right)-t \beta(t+\beta(\lambda-\alpha))\|\nabla F(x(t))\|^{2}\\
&\leqslant\frac{\gamma \lambda-2}{2 t}\left\|\lambda\left(x(t)-x^{*}\right)\right\|^{2} \\
&+\left(\frac{\lambda+1-\alpha}{t}+\frac{\gamma \lambda-2}{2 t}\right)\|t(\dot{x}(t)+\beta \nabla F(x(t)))\|^{2} \\
&+\left(\frac{\lambda+1-\alpha}{t}+\frac{\gamma \lambda-2}{t}\right)\left\langle \lambda(x(t)-x^{*}), t(\dot{x}(t)+\beta \nabla F(x(t)))\right\rangle \\
&+\beta(\alpha-\lambda)\left(F(x(t))-F^{*}\right)-t \beta(t+\beta(\lambda-\alpha))\|\nabla F(x(t))\|^{2}\end{aligned}$$
Noticing that $2(\lambda-\alpha)+\gamma\lambda=0$, we get that:
$$\begin{aligned}
\mathcal{E}^{\prime}(t)+\frac{\gamma \lambda-2}{t} \mathcal{E}(t)&\leqslant \frac{\gamma \lambda-2}{2 t}\left\|\lambda\left(x(t)-x^{*}\right)\right\|^{2} \\
&+\frac{\lambda+\gamma \lambda-\alpha-1}{t}\left\langle\lambda\left(x(t)-x^{*}\right), t(\dot{x}(t)+\beta \nabla F(x(t)))\right\rangle\\
&+\beta(\alpha-\lambda)\left(F(x(t))-F^{*}\right)-t \beta(t+\beta(\lambda-\alpha))\|\nabla F(x(t))\|^{2}.
\end{aligned}
$$
Consequently,
\begin{equation}
\begin{aligned}
\mathcal{E}^{\prime}(t)+\frac{\gamma \lambda-2}{t} \mathcal{E}(t) & \leqslant K(\alpha)\left( \;\dfrac{\lambda}{t}\left\|x(t)-x^{*}\right\|^{2}+\left\langle x(t)-x^{*}, \dot{x}(t)+\beta \nabla F(x(t))\right\rangle \right)
\\
&+\beta(\alpha-\lambda)\left(F(x(t))-F^{*}\right)-t \beta(t+\beta(\lambda-\alpha))\|\nabla F(x(t))\|^{2},
\end{aligned}
\end{equation}
where $K(\alpha)=\frac{2 \alpha \gamma}{(\gamma+2)^{2}}\left(\alpha-1-\frac{2}{\gamma}\right)$.
\subsection{Proof of Lemma \ref{lem:sharp4}}
\label{subsection:Proof for Lemma 4}
Lemma \ref{lem:sharp3} guarantees that for all $t> \beta(\alpha-\lambda)$:
$$
\begin{aligned}
\mathcal{E}^{\prime}(t)+\frac{\gamma \lambda-2}{t} \mathcal{E}(t) & \leqslant K(\alpha)\left( \;\dfrac{\lambda}{t}\left\|x(t)-x^{*}\right\|^{2}+\left\langle x(t)-x^{*}, \dot{x}(t)+\beta \nabla F(x(t))\right\rangle \right)
\\
&+\beta(\alpha-\lambda)\left(F(x(t))-F^{*}\right)-t \beta(t+\beta(\lambda-\alpha))\|\nabla F(x(t))\|^{2}.
\end{aligned}
$$
By adding $\left(\frac{\gamma \lambda-2}{t+\beta(\lambda-\alpha)}-\frac{\gamma \lambda-2}{t}\right) \mathcal{E}(t)$ to both sides we get:
$$
\begin{aligned}
\mathcal{E}^{\prime}(t)+\frac{\gamma \lambda-2}{t+\beta(\lambda-\alpha)} \mathcal{E}(t) & \leqslant K(\alpha)\left(\frac{\lambda}{t}\left\|x(t)-x^{*}\right\|^{2}+\left\langle x(t)-x^{*}, \dot{x}(t)+\beta \nabla F(x(t))\right\rangle\right) \\
&+\beta(\alpha-\lambda)\left(F(x(t))-F^{*}\right)+\left(\frac{\gamma \lambda-2}{t+\beta(\lambda-\alpha)}-\frac{\gamma \lambda-2}{t}\right) \mathcal{E}(t)\\&-t \beta(t+\beta(\lambda-\alpha))\|\nabla F(x(t))\|^{2}\\
 & \leqslant K(\alpha)\left(\frac{\lambda}{t}\left\|x(t)-x^{*}\right\|^{2}+\left\langle x(t)-x^{*}, \dot{x}(t)+\beta \nabla F(x(t))\right\rangle\right) \\
&+\beta(\alpha-\lambda)\left(F(x(t))-F^{*}\right)+\frac{\beta(\alpha-\lambda)(\gamma\lambda-2)}{t(t+\beta(\lambda-\alpha))} \mathcal{E}(t)\\&-t \beta(t+\beta(\lambda-\alpha))\|\nabla F(x(t))\|^{2}
\end{aligned}
$$
Recall that for all $t> \beta(\alpha-\lambda)$,
$$F(x(t))-F^{*} \leq \dfrac{\mathcal{E}(t)}{t(t+ \beta(\lambda-\alpha))}, $$
and thus:
\begin{equation}
\begin{aligned}
\mathcal{E}^{\prime}(t)+\frac{\gamma \lambda-2}{t+\beta(\lambda-\alpha)} \mathcal{E}(t)& \leqslant K(\alpha)\left(\frac{\lambda}{t}\left\|x(t)-x^{*}\right\|^{2}+\left\langle x(t)-x^{*}, \dot{x}(t)+\beta \nabla F(x(t))\right\rangle\right) \\
&+\left(\frac{\beta(\alpha-\lambda)}{t(t+ \beta(\lambda-\alpha))}+\frac{\beta(\alpha-\lambda)(\gamma\lambda-2)}{t(t+\beta(\lambda-\alpha))}\right) \mathcal{E}(t) \\
&-t \beta(t+\beta(\lambda-\alpha))\|\nabla F(x(t))\|^{2}\\
&\leqslant K(\alpha)\left(\frac{\lambda}{t}\left\|x(t)-x^{*}\right\|^{2}+\left\langle x(t)-x^{*}, \dot{x}(t)+\beta \nabla F(x(t))\right\rangle\right) \\
&+\frac{\beta(\alpha-\lambda)(\gamma \lambda-1)}{t^{2}+t \beta(\lambda-\alpha)} \mathcal{E}(t)-t \beta(t+\beta(\lambda-\alpha))\|\nabla F(x(t))\|^{2}.
\end{aligned}
\label{eq:balise}
\end{equation}
The next step is to find a bound of 
$\frac{\lambda}{t}\left\|x(t)-x^{*}\right\|^{2}+\left\langle x(t)-x^{*}, \dot{x}(t)+\beta \nabla F(x(t))\right\rangle$ depending on $\mathcal{E}(t)$. This will be done by applying the inequalities of the following lemma which is proved in Section \ref{sec:proof_sharp5}.
\begin{lemme}
Let $u\in\R^n$, $v\in\R^n$ and $a>0$. Then,
$$|\langle u,v\rangle|\leqslant\frac{a}{2}\|u\|^2+\frac{1}{2a}\|v\|^2,$$
and $$\|u\|^2\leqslant(1+a)\|u+v\|^2+\left(1+\frac{1}{a}\right)\|v\|^2.$$
\label{lem:sharp5}
\end{lemme}
Lemma \ref{lem:sharp5} ensures that for all $t> \beta(\alpha-\lambda)$ and $\theta>0$,
\begin{equation}
|\langle x(t)-x^{*}, \dot{x}(t)+\beta \nabla F(x(t))\rangle| \leqslant \frac{\sqrt{\mu}}{2}\left\|x(t)-x^{*}\right\|^{2}+\frac{1}{2 \sqrt{\mu}}\|\dot{x}(t)+\beta \nabla F(x(t))\|^{2},
\end{equation}
and 
\begin{align}
t^{2}\|\dot{x}(t)+\beta \nabla F(x(t)))\|^{2} &\leqslant\left(1+\theta \frac{\alpha}{t \sqrt{\mu}}\right)\left\|\lambda\left(x(t)-x^{*}\right)+t (\dot{x}(t)+\beta \nabla F(x(t)))\right\|^{2}\nonumber\\
&+\lambda^{2}\left(1+\frac{t \sqrt{\mu}}{\theta \alpha}\right)\left\|x(t)-x^{*}\right\|^{2}.
\end{align}
Hence, for all $\theta>0$,
$$
\begin{aligned}
&\frac{\lambda}{t}\left\|x(t)-x^{*}\right\|^{2}+\left\langle x(t)-x^{*}, \dot{x}(t)+\beta \nabla F(x(t))\right\rangle \\
&\leqslant \left(\frac{\lambda}{t}+\frac{\sqrt{\mu}}{2}\right)\|x(t)-x^*\|^2+\frac{1}{2\sqrt{\mu}}\|\dot{x}(t)+\beta\nabla F(x(t))\|^2\\
&\leqslant \left(\frac{\lambda^2}{2\sqrt{\mu}t^2}+\frac{\lambda}{t}\left(1+\frac{\lambda}{2\theta\alpha}\right)+\frac{\sqrt{\mu}}{2}\right)\|x(t)-x^*\|^2\\
&+\left(\frac{\theta\alpha}{2 \mu t^3}+\frac{1}{2\sqrt{\mu}t^2}\right) \|\lambda(x(t)-x^*)+t (\dot{x}(t)+\beta \nabla F(x(t)))\|^2\\
&\leqslant \left(\frac{\lambda^2}{\mu^{3/2}t^2}+\frac{2\lambda}{\mu t}\left(1+\frac{\lambda}{2\theta\alpha}\right)+\frac{1}{\sqrt{\mu}}\right)(F(x(t))-F^*)\\
&+\left(\frac{\theta\alpha}{2\mu t^3}+\frac{1}{2\sqrt{\mu}t^2}\right) \|\lambda(x(t)-x^*)+t (\dot{x}(t)+\beta \nabla F(x(t)))\|^2,
\end{aligned}
$$
as $F$ satisfies $\mathcal{G}_\mu^2$ and has a unique minimizer.\\
As $\alpha>\lambda$ we have that $\frac{1}{t}<\frac{1}{t+\beta(\lambda-\alpha)}$ and thus:
$$
\begin{aligned}
&\frac{\lambda}{t}\left\|x(t)-x^{*}\right\|^{2}+\left\langle x(t)-x^{*}, \dot{x}(t)+\beta \nabla F(x(t))\right\rangle\\
&\leqslant\left(\tfrac{\lambda^2}{\mu^{3/2}(t+\beta(\lambda-\alpha))^2}+\tfrac{2\lambda}{\mu (t+\beta(\lambda-\alpha))}\left(1+\tfrac{\lambda}{2\theta\alpha}\right)+\tfrac{1}{\sqrt{\mu}}\right)(F(x(t))-F^*)\\
&+\left(\tfrac{\theta\alpha}{2\mu(t+\beta(\lambda-\alpha))^3}+\tfrac{1}{2\sqrt{\mu}(t+\beta(\lambda-\alpha))^2}\right) \|\lambda(x(t)-x^*)+t (\dot{x}(t)+\beta \nabla F(x(t)))\|^2\\
&\leqslant\left(\tfrac{\lambda^2}{\mu^{3/2}(t+\beta(\lambda-\alpha))^4}+\tfrac{2\lambda}{\mu (t+\beta(\lambda-\alpha))^3}\left(1+\tfrac{\lambda}{2\theta\alpha}\right)+\tfrac{1}{\sqrt{\mu}(t+\beta(\lambda-\alpha))^2}\right)t(t+\beta(\lambda-\alpha))(F(x(t))-F^*)\\
&+\left(\tfrac{\theta\alpha}{\mu(t+\beta(\lambda-\alpha))^3}+\tfrac{1}{\sqrt{\mu}(t+\beta(\lambda-\alpha))^2}\right)\frac{1}{2} \|\lambda(x(t)-x^*)+t (\dot{x}(t)+\beta \nabla F(x(t)))\|^2.
\end{aligned}$$

\noindent The parameter $\theta$ is then defined to ensure that $\tfrac{2\lambda}{\mu (t+\beta(\lambda-\alpha))^3}\left(1+\tfrac{\lambda}{2\theta\alpha}\right)=\tfrac{\theta\alpha}{\mu(t+\beta(\lambda-\alpha))^3}$. This equality is satisfied for $\theta=\frac{2}{\gamma+2}(1+\sqrt{2})$ and this choice leads to the following inequalities:
$$
\begin{aligned}
&\frac{\lambda}{t}\left\|x(t)-x^{*}\right\|^{2}+\left\langle x(t)-x^{*}, \dot{x}(t)+\beta \nabla F(x(t))\right\rangle\\
&\leqslant\tfrac{1}{\mu (t+\beta(\lambda-\alpha))^2}\left(\tfrac{\lambda^2}{\sqrt{\mu}(t+\beta(\lambda-\alpha))^2}+\tfrac{\lambda}{t+\beta(\lambda-\alpha)}\left(1+\sqrt{2}\right)+\sqrt{\mu}\right)t(t+\beta(\lambda-\alpha))(F(x(t))-F^*)\\
&+\tfrac{1}{\mu(t+\beta(\lambda-\alpha))^2}\left(\tfrac{\lambda}{t+\beta(\lambda-\alpha)}(1+\sqrt{2})+\sqrt{\mu}\right) \|\lambda(x(t)-x^*)+t (\dot{x}(t)+\beta \nabla F(x(t)))\|^2\\
&\leqslant \tfrac{1}{\mu (t+\beta(\lambda-\alpha))^2}\left(\tfrac{\lambda^2}{\sqrt{\mu}(t+\beta(\lambda-\alpha))^2}+\tfrac{\lambda}{t+\beta(\lambda-\alpha)}\left(1+\sqrt{2}\right)+\sqrt{\mu}\right)\mathcal{E}(t).
\end{aligned}$$
Coming back to \eqref{eq:balise} we get that,
$$
\begin{aligned}
\mathcal{E}^{\prime}(t)+\frac{\gamma \lambda-2}{t+\beta(\lambda-\alpha)} \mathcal{E}(t)  &\leqslant K(\alpha)\left(\frac{\lambda}{t}\left\|x(t)-x^{*}\right\|^{2}+\left\langle x(t)-x^{*}, \dot{x}(t)+\beta \nabla F(x(t))\right\rangle\right) \\
&+\frac{\beta(\alpha-\lambda)(\gamma \lambda-1)}{t(t+\beta(\lambda-\alpha))} \mathcal{E}(t)-t \beta(t+\beta(\lambda-\alpha))\|\nabla F(x(t))\|^{2}\\
&\leqslant \tfrac{K(\alpha)}{\mu (t+\beta(\lambda-\alpha))^2}\left(\tfrac{\lambda^2}{\sqrt{\mu}(t+\beta(\lambda-\alpha))^2}+\tfrac{\lambda}{t+\beta(\lambda-\alpha)}\left(1+\sqrt{2}\right)+\sqrt{\mu}\right)\mathcal{E}(t)\\
&+\frac{\beta(\alpha-\lambda)(\gamma \lambda-1)}{(t+\beta(\lambda-\alpha))^2} \mathcal{E}(t)-t \beta(t+\beta(\lambda-\alpha))\|\nabla F(x(t))\|^{2}.
\end{aligned}
$$
By defining $C_0=\frac{\beta\sqrt{\mu}(\alpha-\lambda)(\gamma\lambda-1)}{K(\alpha)}$, it can be rewritten:
$$
\begin{aligned}
\mathcal{E}^{\prime}(t)+\tfrac{\gamma \lambda-2}{t+\beta(\lambda-\alpha)} \mathcal{E}(t)&\leqslant \tfrac{K(\alpha)}{\mu (t+\beta(\lambda-\alpha))^2}\left(\tfrac{\lambda^2}{\sqrt{\mu}(t+\beta(\lambda-\alpha))^2}+\tfrac{\lambda}{t+\beta(\lambda-\alpha)}\left(1+\sqrt{2}\right)+\sqrt{\mu}(1+C_0)\right)\mathcal{E}(t)\\
&-t \beta(t+\beta(\lambda-\alpha))\|\nabla F(x(t))\|^2.
\end{aligned}
$$
and finally,
\begin{equation}
\mathcal{E}^{\prime}(t)+\frac{\gamma \lambda-2}{t+\beta(\lambda-\alpha)} \mathcal{E}(t)+t \beta(t+\beta(\lambda-\alpha))\|\nabla F(x(t))\|^2 \leqslant \varphi(t+\beta(\lambda-\alpha)) \mathcal{E}(t)
\end{equation}
where 
$$
\varphi : t\mapsto \frac{K(\alpha)}{\mu t^{2}}\left(\sqrt{\mu}(1+C_0)+\frac{2 \alpha}{(\gamma+2) t}(1+\sqrt{2})+\frac{4 \alpha^{2}}{(\gamma+2)^{2} \sqrt{\mu} t^{2}}\right)
$$

\subsection{Proof of Lemma \ref{lem:H1loc}}

\label{sec:proof_H1loc}

Let $F: \mathbb{R}^{n} \rightarrow \mathbb{R}$ be a convex $C^2$ function with a non empty set of minimizer $X^*$. Let $\delta\in(0,1]$ and $x^*\in X^*$.\\
We introduce the following lemma which is proved in Section \ref{sec:proof_C2}.
\begin{lemme}
\label{lem:C2}
Let $F: \mathbb{R}^{n} \rightarrow \mathbb{R}$ be a $C^2$ function. Then, for all $x\in\R^n$ and $\varepsilon>0$, there exists $\nu>0$ such that for all $y\in B(x,\nu)$:
\begin{equation}
    (1-\varepsilon)(y-x)^TH_F(x)(y-x)\leqslant (y-x)^TH_F(y)(y-x)\leqslant (1+\varepsilon)(y-x)^TH_F(x)(y-x).
\end{equation}
\end{lemme}
As $F$ is a $C^2$ function, Lemma \ref{lem:C2} ensures that there exists $\nu>0$ such that for all $x\in B\left(x^*,\nu\right)$:
\begin{equation}
\left(1-\frac{\delta}{4-\delta}\right)K(x)\leqslant (x-x^*)^T H_F(x) (x-x^*)\leqslant \left(1+\frac{\delta}{4-\delta}\right)K(x),
\label{eq:hess_loc}
\end{equation}
where $K(x)=(x-x^*)^T H_F(x^*) (x-x^*)$.\\
Let $\phi_{x,x^*}$ be defined as follows:
$$
\begin{aligned}
    \phi_{x,x^*}:[0,1]&\to\R\\t&\mapsto F\left(tx+(1-t)x^*\right),
\end{aligned}
$$
for some $x\in B\left(x^*,\nu\right)$. The function $\phi_{x,x^*}$ is twice differentiable and we have that for all $t\in[0,1]$:
\begin{equation*}
    \begin{gathered}
    \phi_{x,x^*}^\prime(t)=(x-x^*)^T\nabla F(tx+(1-t)x^*),\\
    \phi_{x,x^*}^{\prime\prime}(t)=(x-x^*)^TH_F(tx+(1-t)x^*)(x-x^*).
    \end{gathered}
\end{equation*}
By rewriting \eqref{eq:hess_loc} at the point $tx+(1-t)x^*$ for some $t\in[0,1]$ we have:
\begin{equation}
\left(1-\frac{\delta}{4-\delta}\right)\phi_{x,x^*}^{\prime\prime}(0)\leqslant\phi_{x,x^*}^{\prime\prime}(t)\leqslant\left(1+\frac{\delta}{4-\delta}\right)\phi_{x,x^*}^{\prime\prime}(0).\label{eq:deltaH1}
\end{equation}
By integrating the left-hand inequality of \eqref{eq:deltaH1} and noticing that $\phi_{x,x^*}^{\prime}(0)=0$ (since $\nabla F(x^*)=0$), we get that:
$$\forall t\in[0,1],~\left(1-\frac{\delta}{4-\delta}\right)\phi_{x,x^*}^{\prime\prime}(0)t\leqslant\phi_{x,x^*}^{\prime}(t).$$
By integrating the right-hand inequality of \eqref{eq:deltaH1}, we get that:
$$\forall t\in[0,1],~
\phi_{x,x^*}(t)-\phi_{x,x^*}(0)\leqslant \left(1+\frac{\delta}{4-\delta}\right)\phi_{x,x^*}^{\prime\prime}(0)\frac{t^2}{2},
$$
and consequently,
\begin{equation*}
    \forall t\in[0,1],~
\phi_{x,x^*}(t)-\phi_{x,x^*}(0)\leqslant \frac{1}{2-\delta}t\phi_{x,x^*}^{\prime}(t).
\end{equation*}
By choosing $t=1$ and rewriting $\phi_{x,x^*}$ and $\phi_{x,x^*}^\prime$ we deduce that
$$F(x)-F^*\leqslant\frac{1}{2-\delta}\langle\nabla F(x), x-x^*\rangle.$$

\subsection{Proof of Lemma \ref{lem:flat2}}
\label{sec:proof_flat2}

We consider the energy function $\mathcal{E}$ defined for all $t\geqslant t_0$ by:
$$
\begin{aligned}
\mathcal{E}(t) =&\left(t^{2}+t \beta(\lambda-\alpha)\right)\left(F(x(t))-F^{*}\right)+\frac{\xi}{2}\left\|x(t)-x^{*}\right\|^{2}\\
&+\frac{1}{2}\left\|\lambda\left(x(t)-x^{*}\right)+t(\dot{x}(t)+\beta \nabla F(x(t)))\right\|^{2}.
\end{aligned}
$$
Let $v:t\mapsto\lambda\left(x(t)-x^{*}\right)+t(\dot{x}(t)+\beta \nabla F(x(t)))$. The function $v$ is differentiable and we have that:
$$
\begin{aligned}
v^{\prime}(t) &=\lambda \dot{x}(t)+t \ddot{x}(t)+\dot{x}(t)+\beta \nabla F(x(t))+t \beta \nabla^{2} F(x(t)) \dot{x}(t) \\
&=(\lambda+1) \dot{x}(t)+\left(-\alpha \dot{x}(t)-t \beta \nabla^{2} F(x(t)) \dot{x}(t)-t \nabla F(x(t))\right)+\beta \nabla F(x(t))\\
&+t \beta \nabla^{2} F(x(t)) \dot{x}(t)\\
&=(\lambda+1-\alpha) \dot{x}(t)+(\beta-t) \nabla F(x(t)).
\end{aligned}
$$
By differentiating the function $\mathcal{E}(t)$, we get that:
$$
\begin{aligned}
\mathcal{E}^{\prime}(t) =&(2 t+\beta(\lambda-\alpha))\left(F(x(t))-F^{*}\right)+\left(t^{2}+t \beta(\lambda-\alpha)\right)\langle\nabla F(x(t)), \dot{x}(t)\rangle+\left\langle v(t), v^{\prime}(t)\right\rangle\\
&+\xi\left\langle x(t)-x^{*}, \dot{x}(t)\right\rangle.
\end{aligned}
$$
Simple calculations give that:
$$
\begin{aligned}
\left\langle v(t), v^{\prime}(t)\right\rangle &=\left\langle\lambda\left(x(t)-x^{*}\right)+t(\dot{x}(t)+\beta \nabla F(x(t))),(\lambda+1-\alpha) \dot{x}(t)+(\beta-t) \nabla F(x(t))\right\rangle \\
&=\lambda(\lambda+1-\alpha)\left\langle x(t)-x^{*}, \dot{x}(t)\right\rangle+\lambda(\beta-t)\left\langle\nabla F(x(t)), x(t)-x^{*}\right\rangle \\
&+t(\lambda+1-\alpha)\|\dot{x}(t)\|^2+t(\beta-t)\langle\nabla F(x(t)), \dot{x}(t)\rangle \\
&+t \beta(\lambda+1-\alpha)\langle\nabla F(x(t)), \dot{x}(t)\rangle+t \beta(\beta-t)\|\nabla F(x(t))\|^2.
\end{aligned}
$$
Consequently,
$$
\begin{aligned}
\mathcal{E}^{\prime}(t) &=(2 t+\beta(\lambda-\alpha))\left(F(x(t))-F^{*}\right)+\left(t^{2}+t \beta(\lambda-\alpha)\right)\langle\nabla F(x(t)), \dot{x}(t)\rangle \\
&+\lambda(\lambda+1-\alpha)\left\langle x(t)-x^{*}, \dot{x}(t)\right\rangle+\lambda(\beta-t)\left\langle\nabla F(x(t)), x(t)-x^{*}\right\rangle \\
&+t(\lambda+1-\alpha)\|\dot{x}(t)\|^2+t(\beta-t)\langle\nabla F(x(t)), \dot{x}(t)\rangle \\
&+t \beta(\lambda+1-\alpha)\langle\nabla F(x(t)), \dot{x}(t)\rangle+t \beta(\beta-t)\|\nabla F(x(t))\|^2 \\
&+\xi\left\langle x(t)-x^{*}, \dot{x}(t)\right\rangle\\
&=(2 t+\beta(\lambda-\alpha))\left(F(x(t))-F^{*}\right)+2 t \beta(\lambda+1-\alpha)\langle\nabla F(x(t)), \dot{x}(t)\rangle \\
&+\lambda(\lambda+1-\alpha)\left\langle x(t)-x^{*}, \dot{x}(t)\right\rangle+\lambda(\beta-t)\left\langle\nabla F(x(t)), x(t)-x^{*}\right\rangle \\
&+t(\lambda+1-\alpha)\|\dot{x}(t)\|^2+t \beta(\beta-t)\|\nabla F(x(t))\|^2 \\
&+\xi\left\langle x(t)-x^{*}, \dot{x}(t)\right\rangle.
\end{aligned}$$
Then by rearranging the terms we obtain that:
$$\begin{aligned}
\mathcal{E}^\prime(t)&=(2 t+\beta(\lambda-\alpha))\left(F(x(t))-F^{*}\right) \\
&+\frac{\lambda+1-\alpha}{t}\left[t^{2}\|\dot{x}(t)\|^2+2 t^{2} \beta\langle\nabla F(x(t)), \dot{x}(t)\rangle+t^{2} \beta^{2}\|\nabla F(x(t))\|^2\right] \\
&+\frac{\lambda(\lambda+1-\alpha)}{t}\left[\left\langle x(t)-x^{*}, t \dot{x}(t)\right\rangle+\left\langle x(t)-x^{*}, t \beta \nabla F(x(t))\right\rangle\right] \\
&-\beta t(t+\beta(\lambda-\alpha))\|\nabla F(x(t))\|^{2} \\
&-\lambda(t+\beta(\lambda- \alpha))\left\langle\nabla F(x(t)), x(t)-x^{*}\right\rangle \\
&+\xi\left\langle x(t)-x^{*}, \dot{x}(t)\right\rangle\\
&=(2 t+\beta(\lambda-\alpha))\left(F(x(t))-F^{*}\right)+\frac{\lambda+1-\alpha}{t}\|t(\dot{x}(t)+\beta \nabla F(x(t)))\|^{2} \\
&+\frac{\lambda(\lambda+1-\alpha)}{t}\left\langle x(t)-x^{*}, t(\dot{x}(t)+\beta \nabla F(x(t)))\right\rangle \\
&-\beta t(t+\beta(\lambda-\alpha))\|\nabla F(x(t))\|^{2} \\
&-\lambda(t+\beta(\lambda- \alpha))\left\langle\nabla F(x(t)), x(t)-x^{*}\right\rangle\\
&+\xi\left\langle x(t)-x^{*}, \dot{x}(t)\right\rangle\\
&=(2 t+\beta(\lambda-\alpha))\left(F(x(t))-F^{*}\right)+\frac{\lambda+1-\alpha}{t}\|\lambda(x(t)-x^{*})+t(\dot{x}(t)+\beta \nabla F(x(t)))\|^{2} \\
&-\frac{\lambda^2(\lambda+1-\alpha)}{t}\|x(t)-x^{*}\|^{2}-\beta t(t+\beta(\lambda-\alpha))\|\nabla F(x(t))\|^{2}  \\
&-\lambda(t+\beta(2(\lambda-\alpha)+1))\left\langle\nabla F(x(t)), x(t)-x^{*}\right\rangle+(\xi-\lambda(\lambda+1-\alpha))\left\langle x(t)-x^{*}, \dot{x}(t)\right\rangle.
\end{aligned}
$$
Given this expression, we can write that:
$$
\begin{aligned}
\mathcal{E}^{\prime}(t) &=((2-\lambda\gamma_1)t+\beta(\lambda-\alpha-\lambda\gamma_1(2(\lambda-\alpha)+1))))\left(F(x(t))-F^{*}\right)\\
&+\lambda(t+\beta(2(\lambda-\alpha)+1))\left[\gamma_1\left(F(x(t))-F^{*}\right)-\left\langle\nabla F(x(t)), x(t)-x^{*}\right\rangle\right]\\
&+\frac{\lambda+1-\alpha}{t}\|\lambda(x(t)-x^{*})+t(\dot{x}(t)+\beta \nabla F(x(t)))\|^{2}-\frac{\lambda^2(\lambda+1-\alpha)}{t}\|x(t)-x^{*}\|^{2} \\
&-\beta t(t+\beta(\lambda-\alpha))\|\nabla F(x(t))\|^{2}+(\xi-\lambda(\lambda+1-\alpha))\left\langle x(t)-x^{*}, \dot{x}(t)\right\rangle.
\end{aligned}
$$
As $F$ satisfies the growth condition $\mathcal{H}_{\gamma_1}$, for all $t\geqslant \beta(2(\alpha-\lambda)-1)$,
\begin{equation*}
\lambda(t+\beta(2(\lambda-\alpha)+1))\left[\gamma_1\left(F(x(t))-F^{*}\right)-\left\langle\nabla F(x(t)), x(t)-x^{*}\right\rangle\right] \leqslant 0.
\end{equation*}
Therefore,
$$\begin{aligned}\mathcal{E}^{\prime}(t) &\leqslant((2-\lambda\gamma_1)t+\beta(\lambda-\alpha-\lambda\gamma_1(2(\lambda-\alpha)+1))\left(F(x(t))-F^{*}\right) \\ 
& +\frac{\lambda+1-\alpha}{t}\|\lambda(x(t)-x^{*})+t(\dot{x}(t)+\beta \nabla F(x(t)))\|^{2}\\
& -\frac{\lambda^2(\lambda+1-\alpha)}{t}\|x(t)-x^{*}\|^{2}-\beta t(t+\beta(\lambda-\alpha))\|\nabla F(x(t))\|^{2}.\\
\end{aligned}
$$

\subsection{Proof of Lemma \ref{lem:delta+1}}
\label{sec:proof_delta+1}
Let $F: \mathbb{R}^{n} \rightarrow \mathbb{R}$ be a convex function having a non empty set of minimizers where $F^*=\inf\limits_{x\in\R^n}F(x)$. Assume that for some $t_1>0$ and $\delta>0$, $F$ satisfies:
\begin{equation}
\int_{t_1}^{+\infty}u^\delta(F(x(u))-F^*)du<+\infty.
\label{eq:ex_delta}
\end{equation}
Let $\varepsilon>0$. Assumption \eqref{eq:ex_delta} ensures that there exists $t_2\geqslant2t_1$ such that:
$$\forall t\geqslant t_2,\quad \int_{t/2}^{t}u^\delta(F(x(u))-F^*)du<\varepsilon.$$
Let $z$ be defined as follows:
$$z:t\mapsto\frac{\int_{t/2}^tu^\delta x(u)du}{\int_{t/2}^tu^\delta du}.$$
Let $t\geqslant t_2$. We define $\nu$ as :
$$ \begin{aligned}
\nu : \mathcal{B}([t/2,t])&\rightarrow[0,1]\\
A\qquad&\mapsto\frac{\int_Au^\delta du}{\int_{t/2}^tu^\delta du},
\end{aligned}$$
where $\mathcal{B}([t/2,t])$ is the Borel $\sigma$-algebra on $[t/2,t]$. Then, we can write that $z(t)=\int_{t/2}^t x(u)d\nu(u)$. As $\nu([t/2,t])=1$ and $F$ is a convex function, Jensen's inequality ensures that:
$$\begin{aligned}
F(z(t))-F^*&=F\left(\int_{t/2}^t x(u)d\nu(u)\right)-F^*\\
&\leqslant \int_{t/2}^t F(x(u))d\nu(u)-F^*\\
&\leqslant \int_{t/2}^t \left(F(x(u))-F^*\right)d\nu(u)\\
&\leqslant \frac{\varepsilon}{\int_{t/2}^tu^\delta du}\\
\end{aligned}$$
Hence, as $t$ tends towards $+\infty$, $F(z(t))-F^*=o\left(t^{-\delta-1}\right).$

\subsection{Proof of Lemma \ref{lem:delta+1phi}}
\label{subsection:delta+1phi}
Let $\phi: \mathbb{R}^{n} \rightarrow \mathbb{R}^+$ such that for some $t_1>0$ and $\delta>0$, $\phi$ satisfies:
\begin{equation}
\int_{t_1}^{+\infty}u^\delta\phi(x(u))du<+\infty.
\label{eq:delta+12}
\end{equation}
Let $\varepsilon>0$. Assumption \eqref{eq:delta+12} guarantees that there exists $t_2\geqslant 2t_1$ such that
\begin{equation*}
\forall t\geqslant t_2,\quad \int_{t/2}^t u^\delta\phi(x(u))du<\varepsilon.
\end{equation*}
Consequently, for all $t\geqslant t_2$,
\begin{equation*}
\inf\limits_{u\in[t/2,t]}\phi(x(u))\int_{t/2}^t u^\delta du<\varepsilon,
\end{equation*}
and
\begin{equation*}
\inf\limits_{u\in[t/2,t]}\phi(x(u))<\frac{\varepsilon(\delta+1)}{t^{\delta+1}-\left(\frac{t}{2}\right)^{\delta+1}}.
\end{equation*}
Hence, as $t\rightarrow+\infty$,
\begin{equation}
\inf\limits_{u\in[t/2,t]}\phi(x(u))=o\left(t^{-\delta-1}\right).
\end{equation}
We recall that $\liminf\limits_{t\rightarrow+\infty} f(t)=\lim\limits_{t\rightarrow+\infty}\left[\inf\limits_{\tau\geqslant t}f(\tau)\right]$. As $\phi$ is a positive function, we get that:
\begin{equation*}
\liminf\limits_{t\rightarrow+\infty} t^{\delta+1}\log(t)\phi(x(t))=l\geqslant0.
\end{equation*}
Suppose that $l>0$. Then there exists $\hat t>t_1$ such that:
\begin{equation*}
\forall t\geqslant \hat t,\quad t^{\delta+1}\log(t)\phi(x(t))\geqslant \frac{l}{2},
\end{equation*}
and hence:
\begin{equation*}
\forall t\geqslant \hat t,\quad t^\delta\phi(x(t))\geqslant \frac{l}{2t\log(t)}.
\end{equation*}
This inequality can not hold as we assume that \eqref{eq:delta+12} is satisfied. We can deduce that $l=0$.

\subsection{Proof of Lemma \ref{lem:sharp5}}
\label{sec:proof_sharp5}

Let $u\in\R^n$, $v\in\R^n$ and $a>0$. The first inequality comes from the following inequalities:
\begin{equation*}
    \begin{aligned}
        \langle u,v\rangle&=\frac{1}{2}\left\|\sqrt{a}u-\frac{v}{\sqrt{a}}\right\|^2-\frac{a}{2}\|u\|^2-\frac{1}{2a}\|v\|^2\geqslant -\frac{a}{2}\|u\|^2-\frac{1}{2a}\|v\|^2,
    \end{aligned}
\end{equation*}
and
\begin{equation*}
    \begin{aligned}
        \langle u,v\rangle&=\frac{a}{2}\|u\|^2+\frac{1}{2a}\|v\|^2-\frac{1}{2}\left\|\sqrt{a}u+\frac{v}{\sqrt{a}}\right\|^2\leqslant \frac{a}{2}\|u\|^2+\frac{1}{2a}\|v\|^2.
    \end{aligned}
\end{equation*}
The second inequality is proved by rewriting $\|u\|^2$ as follows:
\begin{equation*}
    \|u\|^2=\|u+v\|^2+\|v\|^2-2\langle u+v,v\rangle,
\end{equation*}
and by applying the first inequality to $\langle u+v,v\rangle$.

\subsection{Proof of Lemma \ref{lem:C2}}
\label{sec:proof_C2}

Let $F: \mathbb{R}^{n} \rightarrow \mathbb{R}$ be a $C^2$ function. We denote the second order partial derivatives of $F$ by $\partial_{ij} F =\frac{\partial^2 F}{\partial x_i\partial x_j}$ for all $(i,j)\in\llbracket 1,n\rrbracket^2$. 

Let $x\in\R^n$ and $\varepsilon>0$. For all $(i,j)\in\llbracket 1,n\rrbracket^2$, $\partial_{ij} F$ is continuous on $\R^n$ and consequently,
\begin{equation*}
    \exists \tilde\nu>0,~ \forall y\in B(x,\tilde\nu),~(1-\varepsilon)\partial_{ij} F(x)\leqslant \partial_{ij} F(y)\leqslant (1+\varepsilon)\partial_{ij} F(x).
\end{equation*}
By taking the minimal value of $\tilde\nu$ for all $(i,j)\in\llbracket 1,n\rrbracket^2$, we get that there exists $\tilde\nu>0$ such that:
\begin{equation}
    \forall (i,j)\in\llbracket 1,n\rrbracket^2,~\forall y\in B(x,\tilde\nu),~(1-\varepsilon)\partial_{ij} F(x)\leqslant \partial_{ij} F(y)\leqslant (1+\varepsilon)\partial_{ij} F(x).\label{eq:cont_12}
\end{equation}
Let $\nu=\min\left\{\tilde\nu,\left(n\max\limits_{(i,j)\in\llbracket 1,n\rrbracket^2}|\partial_{ij}F(x)|\right)^{-\frac{1}{2}}\right\}$, $y\in B(x,\nu)$ and $h=y-x$. Equation \eqref{eq:cont_12} gives us that for all $(i,j)\in\llbracket 1,n\rrbracket^2$:
\begin{equation}
    \partial_{ij} F(x)h_ih_j-\varepsilon |\partial_{ij} F(x)h_ih_j|\leqslant \partial_{ij} F(y)h_ih_j\leqslant \partial_{ij} F(x)h_ih_j+\varepsilon |\partial_{ij} F(x)h_ih_j|.\label{eq:cont_13}
\end{equation}
We recall that for all $(i,j)\in\llbracket 1,n\rrbracket^2$, $\left(H_F(x)\right)_{i,j}=\partial_{ij}F(x)$ and therefore: $$\forall (x,h)\in\R^n\times\R^n,~h^T H_F(x) h=\sum_{i=1}^n\sum_{j=1}^n\partial_{ij}F(x)h_ih_j.$$
By summing \eqref{eq:cont_13} for all $(i,j)\in\llbracket 1,n\rrbracket^2$, we get that:
\begin{equation*}
    h^TH_F(x)h-\varepsilon\sum_{i=1}^n\sum_{j=1}^n|\partial_{ij} F(x)h_ih_j|\leqslant h^TH_F(y)h\leqslant h^TH_F(x)h+\varepsilon\sum_{i=1}^n\sum_{j=1}^n|\partial_{ij} F(x)h_ih_j|.
\end{equation*}
Noticing that $|h_ih_j|\leqslant \frac{1}{2}\left(h_i^2+h_j^2\right)$ for all $(i,j)\in\llbracket 1,n\rrbracket^2$, we can deduce that:
\begin{equation*}
    \begin{aligned}
    \sum_{i=1}^n\sum_{j=1}^n|\partial_{ij} F(x)h_ih_j|&\leqslant \max\limits_{(i,j)\in\llbracket 1,n\rrbracket^2}|\partial_{ij}F(x)|\sum_{i=1}^n\sum_{j=1}^n|h_ih_j|\\&\leqslant n\max\limits_{(i,j)\in\llbracket 1,n\rrbracket^2}|\partial_{ij}F(x)|\|h\|^2\\&\leqslant n\max\limits_{(i,j)\in\llbracket 1,n\rrbracket^2}|\partial_{ij}F(x)|\nu^2\\&\leqslant 1.
    \end{aligned}
\end{equation*}
Hence, 
\begin{equation*}
    (1-\varepsilon) h^TH_F(x)h\leqslant h^TH_F(y)h\leqslant (1+\varepsilon)h^TH_F(x)h.
\end{equation*}

\section*{Acknowledgements}
The authors acknowledge the support of the French Agence Nationale de la Recherche (ANR) under reference ANR- PRC-CE23 MaSDOL and the support of FMJH Program PGMO 2019-0024 and from the support to this program from EDF-Thales-Orange.

\bibliographystyle{abbrv}
\bibliography{ref.bib}

\end{document}